\documentclass[12pt]{article}
\usepackage{amssymb}
\usepackage{amsfonts}
\usepackage{amsmath}
\usepackage{amsmath}
\usepackage{verbatim}
\usepackage{color}
\usepackage{showlabels}
\usepackage[active]{srcltx}
\include{srctex}
\allowdisplaybreaks[4]
\newtheorem{theorem}{Theorem}[section]

\newtheorem{definition}[theorem]{Definition}

\newtheorem{problem}[theorem]{Problem}

\newtheorem{remark}[theorem]{Remark}

\let\Section=\section
\def\section{\setcounter{equation}{0}\Section}
\newenvironment{proof}[1][Proof]{\textbf{#1.} }{\ \rule{0.5em}{0.5em}}
\newcommand{\R}{\mathbb{R}}

\def\RR{\mathbb{R}}

\def\EE{\mathbb{E}}

\def\si{{\sigma}}

\def\si{{\sigma}}

\def\vare{{\varepsilon}}
\def \eref#1{\hbox{(\ref{#1})}}

\begin{document}

\title{Smoothness of the joint density for spatially homogeneous SPDEs
\thanks{Y. Hu is partially supported by a grant from the Simons Foundation \# 209206. D. Nualart is supported by the
NSF grant  DMS 1208625. X. Sun is supported by the NNSF of China (No.: 11271169) and the Priority Academic Program Development of Jiangsu Higher Education Institutions.\newline
Keywords and phrases: Stochastic partial differential equations, Malliavin calculus,
spatially homogeneous covariances, smoothness of joint density, strict positivity. }}

\author{{Yaozhong Hu}$^a$\footnote{E-mail: yhu@ku.edu}\;
\  \ \ {Jingyu Huang}$^a$\footnote{E-mail: huangjy@ku.edu}\;
\  \ \ {David Nualart}$^a$\footnote{E-mail: nualart@ku.edu}\;
\  \ \ {Xiaobin Sun}$^b$\footnote{E-mail: xbsun@mail.nankai.edu.cn}
\\
 \small  a. Department of Mathematics, University of Kansas, Lawrence, Kansas, 66045, USA.\\
 \small  b. School of Mathematical Sciences, Nankai University, Tianjin 300071, China.}\,

\date{}

\maketitle

\begin{abstract}
In this paper we consider a general class of second order stochastic partial
differential equations on $\mathbb{R}^d$ driven by a Gaussian noise which is
white in time and it has a homogeneous spatial covariance. Using the techniques of
Malliavin calculus we derive the smoothness of the density of the solution at a fixed
number of points  $(t,x_1), \dots, (t,x_n)$, $t>0$, assuming some suitable regularity and
non degeneracy assumptions. We also prove that the density is strictly positive in the
interior of the support of the law.
\end{abstract}

\section{Introduction}
Consider the stochastic partial differential equation
\begin{equation}\label{Eq}
Lu(t,x)=b(u(t,x))+\sigma(u(t,x))\dot{W}(t,x),
\end{equation}
$t\ge 0$, $x\in \mathbb{R}^d$,
with vanishing initial conditions, where $L$ denotes a second order partial differential operator.
The coefficients $b$ and $\sigma$ are real-valued functions and the noise  $\dot{W}(t,x)$  is a
Gaussian field which is white in time and it has a spatially homogeneous covariance in the space variable.
A mild solution to this equation can be formulated using the
Green kernel  $\Gamma(t,dx)$   associated with the operator $L$ (see Definition  \ref{def mild sol}).
This requires the notion of stochastic integral introduced by Walsh in  \cite{Wal} if  $\Gamma(t,x)$
is a real-valued function or Dalang's extension of Walsh integral
(see \cite{Da}) when $\Gamma$ is a measure.

In  \cite{NQ}, Nualart and Quer-Sardanyons have studied the existence and
smoothness of the density of the solution $u(t,x)$ at a fixed point  $(t,x) \in (0,\infty) \times  \mathbb{R}^d$
using techniques of Malliavin calculus. The smoothness of the density follows from the fact  that
 the norm of the Malliavin derivative of $u(t,x)$ has inverse moments of all orders, assuming some
  suitable non degeneracy and regularity conditions. The basic assumptions  are that $b$ and $\sigma$
  are smooth  with bounded partial derivatives of all orders, $|\sigma(z)| \ge c>0$ for all $z$ and
\begin{equation} \label{Eq2}
C \varepsilon^\eta \le \int_0^\varepsilon \int_{\mathbb{R}^d} |\mathcal{F} \Gamma(r) (\xi) |^2 \mu(d\xi) dr <\infty
\end{equation}
for some $\eta>0$ and all $\varepsilon$ small enough, where $\mu$ is the spectral measure of the noise
and $\mathcal{F}$ denotes the Fourier transform.
This general result extends previous work  of Quer-Sardanyons and Sanz-Sol\'e \cite{QS1} for  the case when $L$ corresponds to
 the three dimensional wave equation.

The purpose of this paper is to establish the
smoothness of the joint density of the solution to equation (\ref{Eq}) at a fixed number of points
$(t,x_1), \dots, (t,x_n)$, where $t>0$ and $x_i \in \mathbb{R}^d$. This kind of problem was studied by
Bally and Pardoux  in \cite{BP} for the one-dimensional stochastic heat equation driven by a
space-time white noise.   The extension of  this result to equation (\ref{Eq}) presents  new difficulties
and requires additional  non degeneracy conditions,  in addition to (\ref{Eq2}), because we need to
handle the determinant of the Malliavin matrix of the random vector  $u(t,x_1), \dots,  u(t,x_n)$.
The basic ingredient is to impose that leading terms as $\vare\rightarrow 0$ in  the matrix
\[
 \left(  \int_0^\varepsilon \int_{\mathbb{R}^d}  \langle  \Gamma(r, *+x_j)  ,\Gamma(r, *+x_i  )
 \rangle_{\mathcal{H}}  dr   \right)_{ 1\le i,j \le n }
\]  
are the diagonal ones  given by  \eref{Eq2}
  (see the hypotheses  ({\bf H3}) and  ({\bf H4}) below).    These hypotheses are
 related, although different, to the ones imposed by Nualart in \cite{NuaE} to establish
 the smoothness of the density for the solution of a system of SDPEs.

The paper is organized as follows. After some preliminaries, Section 3 is devoted to the proof of the smoothness of the density of the vector $u(t,x_1), \dots,  u(t,x_n)$. In Section 4 we derive the positivity of the density in the interior of the support following the general criterion established by  Nualart in \cite{NuaD}.  Finally, in Section 5 we apply these results to the basic examples of the stochastic heat and wave equations and to the spatial covariances given by the Riesz, Bessel and fractional kernels.

\section{Preliminaries}\label{preliminaries}

Consider a  non-negative and non-negative definite function $f$ which is continuous on $\RR^d \setminus \{0\}$.
We  assume that $f$ is the Fourier transform
of a non-negative tempered measure $\mu$ on $\RR^d$ (called the
spectral measure of $f$). That is, for all $\varphi$ belonging to
the space $\mathcal{S}(\RR^d)$ of rapidly decreasing $\mathcal{C}^{\infty}$
functions on $\RR^d$
\begin{equation}\label{def of spectral measure mu}
\int_{\RR^d}f(x)\varphi(x)dx=\int_{\RR^d}\mathcal{F}\varphi(\xi)\mu(d\xi),
\end{equation}
and there is an integer $m\geq 1$ such that
\begin{equation}
\int_{\RR^d}(1+|\xi|^2)^{-m}\mu(d\xi)<\infty\,.\nonumber
\end{equation}
Here we have denoted by $\mathcal{F}\varphi$ the Fourier transform of
$\varphi\in \mathcal{S}(\RR^d)$, given by
$\mathcal{F}\varphi(\xi)=\int_{\RR^d}\varphi(x)e^{- i\xi\cdot x}dx$.

  Let
$\mathcal{C}_0^{\infty}([0,\infty)\times \RR^{d})$ be the space of smooth 
 functions with compact support on $[0,\infty) \times
\R^d$.  Consider a  family of zero mean Gaussian  random variables
$W=\{W(\varphi), \varphi \in \mathcal{C}_0^{\infty}([0,\infty)\times\RR^{d})\}$,
defined on  a complete probability space $(\Omega, \mathcal{F}, \mathbb{P})$,
with covariance
\begin{equation}\label{cov}
\EE(W(\varphi)W(\psi))=\int_0^{\infty}\int_{\RR^d}\int_{\RR^d}\varphi(t,x)f(x-y)\psi(t,y)dxdydt\,.
\end{equation}
The  covariance \eref{cov} can also be written, using Fourier transform, as
\begin{equation}
\EE (W(\varphi)W(\psi))=\int_0^\infty\int_{\RR^d} \mathcal{F}\varphi(t)(\xi) \overline {\mathcal{F}\psi(t)(\xi)}\mu(d\xi)dt\,.\nonumber
\end{equation}

The main assumptions on the differential operator $L$ in \eref{Eq} can be stated as follows:

({\bf H1}) The fundamental solution to $Lu=0$, denoted by $\Gamma$, satisfies that for all $t > 0$ , $\Gamma(t)$ is a nonnegative measure  with  rapid decrease, such that for all $T>0$
\begin{equation}
\int_0^T \int_{\RR^d} |\mathcal{F}\Gamma(t)(\xi)|^2\mu(d\xi)dt <   \infty\,,\nonumber
\end{equation}
and
\begin{equation}
\sup_{t \in [0,T]} \Gamma(t,\mathbb{R}^d)\leq C_T <   \infty\,.\nonumber
\end{equation}

The basic examples we are interested in are the stochastic heat and wave equations. More precisely, it is well-known that if $L$ is the heat operator in $\RR^d$, that is, $L=\frac{\partial }{\partial t}-\frac{1}{2}\Delta$, where $\Delta$ denotes the Laplacian operator in $\RR^d$, or if $L$ is the wave operator in $\RR^d$, $d \in \{1,2,3\}$, i.e., $L = \frac{\partial^2}{\partial t^2}-\Delta$, hypothesis ({\bf H1}) is satisfied if and only if
\begin{equation}
\int_{\RR^d}\frac{\mu(d\xi)}{1+|\xi|^2} < \infty\,.\nonumber
\end{equation}

Let $\mathcal{H}$ be the Hilbert space obtained by the completion of
$\mathcal{C}_0^{\infty}(\mathbb{R}^d)$ endowed with the inner product
\begin{equation} \label{def H}
\langle\varphi,\psi\rangle_{\mathcal{H}}=\int_{\RR^d}dx\int_{\RR^d}dy\varphi(x)f(x-y)\psi(y)=\int_{\RR^d} \mathcal{F}(\varphi)(\xi) \overline {\mathcal{F}(\psi)(\xi)}\mu(d\xi),
\end{equation}
$\varphi,\psi \in \mathcal{C}_0^{\infty}(\RR^d)$. Notice that $\mathcal{H}$ may contain distributions.  Set $\mathcal{H}_0=L^2([0,\infty);\mathcal{H})$.

Walsh's classical theory of stochastic integration developed in
\cite{Wal} cannot be applied directly to the mild formulation of
equation \eref{Eq} since $\Gamma$ may not be absolutely continuous with respect to the Lebesgue measure.  We shall use the stochastic integral defined in    \cite[Section 2.3]{DQ} (see also  \cite[Section 3]{NQ}).   We briefly review  the construction and properties
of this integral.

The Gaussian family $W$ can be extended to the space $\mathcal{H}_0$ and we  denote by
$W(g)$ the Gaussian random variable associated with an element $g
\in \mathcal{H}_0$.  It is obvious that ${\bf 1}_{[0,t]}h$ is in
$\mathcal{H}_0$ and we set $W_t(h)=W({\bf 1}_{[0,t]}h)$ for any $t\ge 0$ and
$h \in \mathcal{H}$. Then $W=\{ W_t, t\ge 0\}$ is a cylindrical Wiener
process in the Hilbert space $\mathcal{H}$. That is, for any $h \in \mathcal{H}$,
$\{W_t(h), t\ge 0\}$ is a Brownian motion with variance
$t\| h\| ^2_{\mathcal{H}}$, and
\begin{equation*}
\EE(W_t(h)W_s(g))=(s\wedge t)\langle h,g\rangle_{\mathcal{H}} .
\end{equation*}
Let $\mathcal{F}_t$ be the $\sigma$-field generated by the random
variables $\{W_s(h), h\in \mathcal{H}, 0\leq s \leq t\}$ and the $\mathbb{P}$-null
sets. We define the predictable $\sigma$-field as the $\sigma$-field
in $\Omega\times [0,\infty)$ generated by the sets $\{A \times (s,t],
0\leq s< t , A\in \mathcal{F}_s\}$. Then we can define the
stochastic integral of  an $\mathcal{H}$-valued square-integrable predictable
process $g \in L^2(\Omega\times[0,\infty);\mathcal{H})$ with respect to the
cylindrical Wiener process $W$, denoted by
\begin{equation*}
g\cdot W=\int_0^\infty\int_{\RR^d}g(t,x)W(dt,dx),
\end{equation*}
and we have the isometry property
\begin{equation}\label{isometry property}
\EE|g\cdot W|^2=\EE\int_0^\infty \| g(t)\|^2_{\mathcal{H}} dt.
\end{equation}

Using the above notion of stochastic integral one can introduce the
following definition:

\begin{definition}\label{def mild sol}
A real-valued predictable stochastic process $u=\{u(t, x), t\ge 0
 , x\in \RR^d\}$ is a mild   solution of equation \eref{Eq} if for all $t\ge 0$, $x\in \RR^d$,
\begin{eqnarray*}
 u(t, x)
&=&\int_0^t\int_{\RR^d}\Gamma(t-s,x-y)\si(u(s,y)) W(ds,dy)\\
&&+\int_0^t\int_{\RR^d}b(u(s,x-y))\Gamma(t-s,dy)ds \,\ \ \ \
a.s.
\end{eqnarray*}
\end{definition}

Now we state the existence and uniqueness result of the solution to equation \eref{Eq}. For a proof of this result, see, for example,   \cite[Theorem 4.3]{DQ}.

\begin{theorem}\label{existence and uniqueness THM}
Suppose hypothesis ({\bf H1}) holds, and $\sigma$, $b$ are Lipschitz continuous.
Then there exists a unique mild    solution $u$ to equation
\eref{Eq} such that for all $p\geq 1$ and $T>0$,
\begin{equation}
\sup_{(t,x) \in [0,T]\times\mathbb{R}^d}\EE|u(t,x)|^p< \infty.
\label{moments estimate}
\end{equation}
\end{theorem}
 
Next we recall some elements of Malliavin calculus  which will be used to prove the main results of this paper.  We consider the Hilbert space $\mathcal{H}_0$ and the Gaussian family of random variables $\{W(h), h \in \mathcal{H}_0\}$ defined above. Then $\{W(h), h \in \mathcal{H}_0\}$ is a centered Gaussian process such that $\EE (W(h_1)W(h_2))=\langle h_1, h_2 \rangle_{\mathcal{H}_0}, h_1, h_2 \in \mathcal{H}_0$. In this framework we can develop a 
  Malliavin calculus  (see, for instance, \cite{NuaD}). The Malliavin derivative is denoted by $D$ and for any $N\geq 1$ and any real number $p\ge 2$, the domain of the iterated derivative $D^N$ in $L^p(\Omega ; \mathcal{H}_0^{\otimes N})$ is denoted by $\mathbb{D}^{N, p}$. We shall also use the notation
\begin{equation*}
\mathbb{D}^{\infty}=\cap_{p \geq 1} \cap_{k \geq 1} \mathbb{D}^{k,p}\,.
\end{equation*}
Note that for any random variable $X$ in the domain of the derivative operator $D$, $DX$ defines an $\mathcal{H}_0$-valued random variable. In particular, for some fixed $r  \ge 0$, $DX(r,*)$ is an element of $\mathcal{H}$, which will be denoted by $D_{r,*}X$.

If $x_1, \dots, x_n$ are points in $\RR^d$ we will make use of the notation  $
u(t,\underline{x})=(u(t,x_1),\ldots, u(t,x_n))$.
In order to study the   smoothness and strict positivity of the (joint) 
density of a random vector of the form $u(t,\underline{x})$, we need to assume some moment estimates for the increments of the solution. We will also need to assume some integral bounds of the fundamental solution $\Gamma$. We list these assumptions below.

\medskip
({\bf H2})  There exist positive constants $\kappa_1$ and $ \kappa_2$ such that for all $s,t  \in [0,T]$, $x,y \in \RR^d$, $T>0$ and $p \geq 1$,
\begin{eqnarray}
&&\EE |u(s, x)-u(t,x)|^p \leq C_{p,T} |t-s|^{\kappa_1 p}\,, \label{Holder time}\\
&& \EE |u(t,x)-u(t,y)|^p \leq C_{p, T}|x-y|^{\kappa_2 p} \label{Holder space}
\end{eqnarray}
for some constant $C_{p,T}$ which only depends on $p,T$.

\medskip
({\bf H3}) There exists $\eta >0$ and $\varepsilon_0>0$ such that for all $0<\varepsilon\le \varepsilon_0$,
\begin{equation}
C \varepsilon^{\eta} \leq \int_0^{\varepsilon}  \|\Gamma(r)\|^2_{\mathcal{H}} dr\nonumber
\end{equation}
for some constant $C> 0$.

\medskip
({\bf H4}) Let $\eta$ be given in hypothesis ({\bf H3}) and  $\kappa_1$ and $\kappa_2$ be given in ({\bf H2}).

 (i) There exists $\eta_1  > \eta$ and $\varepsilon_1>0$ such that for all $0<\varepsilon\le \varepsilon_1$,
\begin{equation}\label{hyp 4-2}
\int_0^{\varepsilon} r^{\kappa_1} \|\Gamma(r)\|^2_{\mathcal{H}}dr\leq C\varepsilon^{\eta_1}.
\end{equation}

(ii)  There exists $\eta_2  > \eta$ such that for each fixed non zero $w\in \RR^d$,  there exists a positive constant $C_w$ and $\varepsilon_2>0$     satisfying
\begin{equation}\label{hyp 4-1}
\int_0^{\varepsilon} \langle \Gamma(r,*), \Gamma(r,w+*)\rangle_{\mathcal{H}}dr \leq C_w \varepsilon^{\eta_2}\,
\end{equation}
for all $0<\varepsilon\le \varepsilon_2$.

(iii)
The  measure $\Psi(t)$ defined by  $|x|^{\kappa_2} \Gamma(t,dx)$ satisfies $\int_0^T \int_{\RR^d} | \mathcal{F} \Psi(t) (\xi) |^2 \mu(d\xi) dt<\infty$  and
there exists $\eta_3>\eta$ such that for each fixed $w\in \RR^d$,
there exists a positive constant $C_w$ and $\varepsilon_3>0$   satisfying
\begin{equation}\label{hyp 4-3}
\int_0^{\varepsilon} \langle  |*|^{\kappa_2} \Gamma(r,*), \Gamma(r, w+*)\rangle_{\mathcal{H}}dr\leq C_{w}\varepsilon^{\eta_3}
\end{equation}
 for all $0<\varepsilon\le \varepsilon_3$.

Along the paper, $C_p$ and $C$ will denote generic constants which may change from line to line and $C_p$ depends on $p\geq 2$.

\section{Existence and smoothness of the density}

Fix $t>0$ and fix   distinct points $x_1,\ldots, x_n$ of $\RR^d$. Let $u(t,x)$  denote the solution of equation (\ref{Eq}).
Recall that  $
u(t,\underline{x})=(u(t,x_1),\ldots, u(t,x_n))$.
In this section we give sufficient conditions for the existence and smoothness of the density of the law of the random vector $u(t,\underline{x})$,
using    Malliavin calculus. The main result is the following theorem.
\begin{theorem} \label{THM smooth density}
Assume that conditions ({\bf H1})-({\bf H4}) hold, and the coefficients $\sigma$, $b$ are $\mathcal{C}^{\infty}$ functions with bounded derivatives of
all orders. Assume that  there exists a positive  constant $C_1$ such that $|\sigma(u(t,x_i))| \ge C_1$ $\mathbb{P}$-a.s. for any $i=1,\dots,n$.
Then the law of the random vector $u(t,\underline{x})$ has a $\mathcal{C}^{\infty}$ density with respect to the Lebesgue measure on $\RR^n$.
\end{theorem}

\begin{remark}
Using a localization procedure developed in   \cite[Theorem 3.1]{BP}, we can prove a version of Theorem \ref{THM smooth density} without assuming that
$|\sigma(u(t, x_i))|\geq C_1$ $\mathbb{P}$-a.s., for any $i=1,\ldots, n$. In this case, we conclude that the law of $u(t,\underline {x})$ has a smooth density on
$\{y\in \RR: \sigma(y) \neq 0\}^n$.
\end{remark}
\begin{proof}
We begin by noting that according to Proposition 6.1 in \cite{NQ}, for each fixed $(t,x) \in [0,\infty)\times \RR^d$, $u(t,x) \in \mathbb{D}^{\infty}$. If we denote by  $M_t(\underline{x})$ the Malliavin covariance matrix \newline $(\langle Du(t,x_i), Du(t,x_j)\rangle_{\mathcal{H}_0})_{1\leq i,j\leq n}$, then,
taking in to account  Theorem 2.1.4 in \cite{NuaD}, we only need to show that the determinant of the Malliavin covariance matrix of $u(t,\underline{x})$ has negative moments of all orders, that is
$$
\EE \left( \text {det} M_t(\underline{x})\right)^{-p}< \infty
$$
for all $p\geq 2$. It suffices to check that for any $p\geq 2$, there exists an $ \delta_0(p)>0$ such that for all $0 < \delta\leq \delta_0(p)$
$$
\mathbb{P}\{\text{det} M_t(\underline{x})\leq \delta\}\leq C\delta^{  p},
$$
for some   constant $C$ not depending on $\delta$.

We begin by noting that
\begin{equation}
\text{det} M_t(\underline{x})\geq \left(\inf_{  \|\xi\|=1} \xi^{T}M_t(\underline{x})\xi  \right)^n.  \label{Mt}
\end{equation}
The derivative of the solution satisfies the following equation (see Proposition 5.1 in \cite{NQ})
\begin{eqnarray*}
D_{r, \ast}u(t, x)&=&\Gamma(t-r, x-\ast)\sigma\left(u(r,\ast)\right)+
\int_r^t\int_{\RR^d}\Gamma(t-s,x-y)\sigma^{\prime}(u(s,y))D_{r, \ast}u(s, y)W(ds,dy)\\
&&+\int_r^t\int_{\RR^d}b^{\prime}(u(s,x-y))D_{r, \ast}u(s, x-y)\Gamma(t-s,dy)ds.
\end{eqnarray*}
Therefore, we can write
$$
\xi^{T}M_t(\underline{x})\xi\geq \int^t_{t-\varepsilon}\|\sum^n_{i=1}D_{r, \ast}u(t, x_i)\xi_i\|^2_{\mathcal{H}}dr
\geq  \frac{1}{2}\mathcal{A}_1- \mathcal{A}_2,
$$
where
\begin{eqnarray*}
\mathcal{A}_1&=&\int^t_{t-\varepsilon}\|\sum^n_{i=1}\Gamma(t-r, x_i-\ast)\sigma(u(r,\ast))\xi_i\|^2_{\mathcal{H}}dr\,, \\
\mathcal{A}_2&=&\int^t_{t-\varepsilon}\|a(r, t, \underline{x}, \ast)\|^2_{\mathcal{H}}dr\,,
\end{eqnarray*}
and
\begin{eqnarray*}
a(r, t, \underline{x}, \ast)&=&\sum^{n}_{i=1}\int_r^t\int_{\RR^d}\Gamma(t-s,x_i-y)\sigma^{\prime}(u(s,y))D_{r, \ast}u(s, y)W(ds,dy)\xi_i\\
&&+\sum^{n}_{i=1}\int_r^t\int_{\RR^d}b^{\prime}(u(s,x_i-y))D_{r, \ast}u(s, x_i-y)\Gamma(t-s,dy)ds\xi_i\,.
\end{eqnarray*}
The term $\mathcal{A}_1$ can be estimated as follows
\begin{eqnarray*}
\mathcal{A}_1&=&\int^t_{t-\varepsilon}\langle\sum^n_{i=1}\Gamma(t-r, x_i-\ast)\sigma(u(r,\ast))\xi_i,
\sum^n_{j=1}\Gamma(t-r, x_j-\ast)\sigma(u(r,\ast))\xi_j\rangle_{\mathcal{H}}dr\\
&=&\int^t_{t-\varepsilon}\sum^{n}_{i=1}\sum^{n}_{j=1}\xi_i\xi_j\langle \Gamma(t-r, x_i-\ast)\sigma(u(t,x_i)),
\Gamma(t-r, x_j-\ast)\sigma(u(t,x_j))\rangle_{\mathcal{H}}dr\\
&&+\int^t_{t-\varepsilon}\sum^{n}_{i=1}\sum^{n}_{j=1}\xi_i\xi_j\Big[\langle \Gamma(t-r, x_i-\ast)\sigma(u(r,\ast)),
\Gamma(t-r, x_j-\ast)\sigma(u(r,\ast))\rangle_{\mathcal{H}}\\
&&-\langle \Gamma(t-r, x_i-\ast)\sigma(u(t,x_i)),
\Gamma(t-r, x_j-\ast)\sigma(u(t,x_j))\rangle_{\mathcal{H}}\Big]dr\\
&=&\int^t_{t-\varepsilon}\sum^{n}_{i=1}\|\xi_i\|^2\|\Gamma(t-r, x_i-\ast)\sigma(u(t,x_i))\|^2_{\mathcal{H}}dr\\
&&+\int^t_{t-\varepsilon}\sum^{n}_{i,j=1, i\neq j}\xi_i\xi_j\langle \Gamma(t-r, x_i-\ast)\sigma(u(t,x_i)),
\Gamma(t-r, x_j-\ast)\sigma(u(t,x_j))\rangle_{\mathcal{H}}dr\\
&&+\int^t_{t-\varepsilon}\sum^{n}_{i=1}\sum^{n}_{j=1}\xi_i\xi_j\left[\langle \Gamma(t-r, x_i-\ast)\sigma(u(r,\ast)),
\Gamma(t-r, x_j-\ast)\sigma(u(r,\ast))\rangle_{\mathcal{H}}\right.\\
&&-\left.\langle \Gamma(t-r, x_i-\ast)\sigma(u(t,x_i)),
\Gamma(t-r, x_j-\ast)\sigma(u(t,x_j))\rangle_{\mathcal{H}}\right]dr\\
&\geq&\mathcal{A}_{11}-|\mathcal{A}_{12}|-|\mathcal{A}_{13}|\,,
\end{eqnarray*}
where
\begin{eqnarray*}
\mathcal{A}_{11}&=&\int^t_{t-\varepsilon}\sum^{n}_{i=1}\|\xi_i\|^2\|\Gamma(t-r, x_i-\ast)\sigma(u(t,x_i))\|^2_{\mathcal{H}}dr\,,\\
\mathcal{A}_{12}&=&\int^t_{t-\varepsilon}\sum^{n}_{i,j=1, i\neq j}\xi_i\xi_j\langle \Gamma(t-r, x_i-\ast)\sigma(u(t,x_i))\,,
\Gamma(t-r, x_j-\ast)\sigma(u(t,x_j))\rangle_{\mathcal{H}}dr\,,\\
\mathcal{A}_{13}&=&\int^t_{t-\varepsilon}\sum^{n}_{i=1}\sum^{n}_{j=1}\xi_i\xi_j\Big[\langle \Gamma(t-r, x_i-\ast)\sigma(u(r,\ast))\,,
\Gamma(t-r, x_j-\ast)\sigma(u(r,\ast))\rangle_{\mathcal{H}}\\
&&-\langle \Gamma(t-r, x_i-\ast)\sigma(u(t,x_i))\,,
\Gamma(t-r, x_j-\ast)\sigma(u(t,x_j))\rangle_{\mathcal{H}}\Big]dr\,.
\end{eqnarray*}
Then, using the fact that $|\sigma(u(t, x_i))|\geq C_1$, for all $i=1,\ldots, n$, we have
\begin{eqnarray*}
\xi^{T}M_t(\underline{x})\xi&\geq& \frac{1}{2}\mathcal{A}_{11}- \frac{1}{2}|\mathcal{A}_{12}|- \frac{1}{2}|\mathcal{A}_{13}|
- \mathcal{A}_{2}\\
&\geq& \frac{1}{2}C_1 g(\varepsilon)- \frac{1}{2}|\mathcal{A}_{12}|-\frac{1}{2} |\mathcal{A}_{13}|
- \mathcal{A}_{2}\,,
\end{eqnarray*}
where
\begin{equation*}
g(\varepsilon)=\int_0^{\varepsilon}\int_{\RR^d}|\mathcal{F}{\Gamma}(s)(\xi)|^2\mu(d\xi)ds\,.
\end{equation*}
Taking $\varepsilon$ such that $\frac{1}{4}C_1 g(\varepsilon)=\delta^{1/n}$, we obtain
\begin{eqnarray}
&&\mathbb{P}\left\{\inf_{  \|\xi\|=1} \xi^{T}M_t(\underline{x})\xi \leq \delta^{1/n}\right\}\nonumber\\
&\leq&\mathbb{P}\left\{\sup_{  \|\xi\|=1} \left( |\mathcal{A}_{12}| + |\mathcal{A}_{13}|
+  2\mathcal{A}_{2} \right)\geq \frac{1}{2}C_1 g(\varepsilon)\right\}\nonumber\\
&\leq&C_p g(\varepsilon)^{-p}\left[\mathbb{E}\Big(\sup_{  \|\xi\|=1}|\mathcal{A}_{12}|^p\Big)
+\mathbb{E}\Big(\sup_{  \|\xi\|=1}|\mathcal{A}_{13}|^p\Big)
+\mathbb{E}\Big(\sup_{ \|\xi\|=1}|\mathcal{A}_{2}|^p\Big)\right]. \label{123}
\end{eqnarray}
Next, we treat  each of the above expectations   separately. For the first expectation of \eref{123}, using (\ref{moments estimate}) and  property  (ii) in condition ({\bf H4}), we can write
\begin{eqnarray}
&&\mathbb{E}\left(\sup_{  \|\xi\|=1}|\mathcal{A}_{12}|^p\right)\nonumber\\
&=&\mathbb{E}\left(\sup_{  \|\xi\|=1}\left|\int^{\varepsilon}_{0}\sum^{n}_{i,j=1, i\neq j}
\xi_i\xi_j\langle \Gamma(r, x_i-\ast)\sigma(u(t,x_i)),\Gamma(r, x_j-\ast)\sigma(u(t,x_j))\rangle_{\mathcal{H}}dr\right|^p\right)\nonumber\\
&\leq&C_p\sum^{n}_{i,j=1, i\neq j}\left[\mathbb{E}\left(|\sigma(u(t,x_i))\sigma(u(t,x_j))|^p \right)
\left|\int^{\varepsilon}_{0}\langle \Gamma(r, x_i-\ast),\Gamma(r, x_j-\ast)\rangle_{\mathcal{H}}dr\right|^p\right]\nonumber\\
&\leq&C_{p}\varepsilon^{\eta_2 p}. \label{A12}
\end{eqnarray}
For the second expectation of \eref{123},  using Minkowski's inequality and property (i) and (iii) in condition ({\bf H4}), we get
\begin{eqnarray}
&&\mathbb{E}\left(\sup_{  \|\xi\|=1}|\mathcal{A}_{13}|^p\right) \nonumber\\
&\leq&C_p\sum^{n}_{i,j=1}\mathbb{E}\Big|\int^{t}_{t-\varepsilon}dr\int_{\RR^d}\int_{\RR^d}
[\sigma(u(r,z))\sigma(u(r,y))-\sigma(u(t,x_i))\sigma(u(t,x_j))]\nonumber\\
&&\times \Gamma(t-r, x_i-dz)\Gamma(t-r, x_j-dy)f(z-y)\Big|^p\nonumber\\
&\leq&C_p\sum^{n}_{i,j=1}\Big(\int^{t}_{t-\varepsilon}dr\int_{\RR^d}\int_{\RR^d}
\|\sigma(u(r,z))\sigma(u(r,y))-\sigma(u(t,x_i))\sigma(u(t,x_j))\|_{L^{p}(\Omega)}\nonumber\\
&&\times \Gamma(t-r, x_i-dz)\Gamma(t-r, x_j-dy)f(z-y)\Big)^p\nonumber\\
&\leq&C_p\sum^{n}_{i,j=1}\Big(\int^{\varepsilon}_{0}dr\int_{\RR^d}\int_{\RR^d}
(r^{\kappa_1}+|x_i-z|^{\kappa_2}+|x_j-y|^{\kappa_2})\nonumber\\
&&\!\!\times\Gamma(r, x_i-dz)\Gamma(r, x_j-dy)f(z-y)\Big)^p\nonumber\\
&\leq&C_{p}\left|\int^{\varepsilon}_{0}r^{\kappa_1}\|\Gamma(r, \ast)\|^2_{\mathcal{H}}dr\right|^p+C_p\sum^{n}_{i,j=1}\left|\int^{\varepsilon}_{0}\langle|*|^{\kappa_2}\Gamma(r, \ast), \Gamma(r, x_j-x_i+\ast)\rangle_{\mathcal{H}}dr\right|^p\nonumber\\
&\leq&C_{p}\varepsilon^{\eta_1 p}+C_{p}\varepsilon^{\eta_3 p}. \label{A13}
\end{eqnarray}
Finally, we treat the last expectation of \eref{123} and we obtain the following inequalities
\begin{eqnarray*}
&&\mathbb{E}\left(\sup_{ \|\xi\|=1}|\mathcal{A}_{2}|^p\right)\\
&& \leq
C_p\mathbb{E}\left(\sup_{ \|\xi\|=1}
\int^{t}_{t-\varepsilon}\left\|\sum^{n}_{i=1}\int_r^t\int_{\RR^d}\Gamma(t-s,x_i-y)\sigma^{\prime}(u(s,y))D_{r, \ast}u(s, y)W(ds,dy)\xi_i\right\|^2_{\mathcal{H}}dr\right)^p\\
&&+C_p\mathbb{E}\left(\sup_{ \|\xi\|=1}
\int^{t}_{t-\varepsilon}\left\|\sum^{n}_{i=1}\int_r^t\int_{\RR^d}b^{\prime}(u(s,x_i-y))D_{r, \ast}u(s, x-y)\Gamma(t-s,dy)ds\xi_i\right\|^2_{\mathcal{H}}dr\right)^p\\
& & := \quad T_1+T_2.
\end{eqnarray*}
For any $\varphi$, $\psi$ in $\mathcal{H}_0$ we use the notation
\begin{equation*}
\langle \varphi, \psi \rangle_{\mathcal{H}_{t-\varepsilon,t}}:= \int_{t-\varepsilon}^t \langle \varphi(s, \ast), \psi(s, *)\rangle_{\mathcal{H}} ds.
\end{equation*}
Using equation (3.13) and the inequality (5.26)  in \cite{NQ}, we obtain
\begin{eqnarray}
T_1&\leq& C_p \sum_{i=1}^n \EE \left \| \int _{t-\varepsilon} ^t \int_{\RR^d} \Gamma(t-s,x_i-y) \sigma^{\prime} (u (s,x_i-y)) Du(s,x-y)W(ds,dy)\right\|^{2p}_{\mathcal{H}_{t-\varepsilon,t}} \nonumber\\
&\leq& g(\varepsilon)^p \sup_{t-\varepsilon \leq s \leq t, x \in \RR^d} \EE \left\|D u(s,x) \right\|^{2p}_{\mathcal{H}_{t-\varepsilon,t}} \nonumber\\
&\leq&C_{p}g(\varepsilon)^{2p} \label{T1}\,.
\end{eqnarray}
 For $T_2$, using Cauchy-Schwartz inequality, our assumption on $b^{\prime}$,  Minkowski's inequality and the estimate (5.26) in \cite{NQ}, we obtain the bound
\begin{eqnarray}
T_2&\leq&C_p\sum^{n}_{i=1}\mathbb{E}\left\|
\int^{t}_{t-\varepsilon}\int_{\RR^d}b^{\prime}(u(s,x_i-y))Du(s, x_i-y)\Gamma(t-s,dy)ds\right\|_{\mathcal{H}_{t-\varepsilon, t}}^{2p}\nonumber\\
&\leq&C_{p}\bigg(\int^{t}_{t-\varepsilon}\int_{\RR^d}\Gamma(t-s,dy)
ds\bigg)^{2p}\sup_{t-\varepsilon\leq s\leq t, x\in\RR^d}\mathbb{E}\|Du(s, x)\|^{2p}_{\mathcal{H}_{t-\varepsilon, t}}\nonumber\\
&\leq&C_{p}g(\varepsilon)^{p}\varepsilon^{2p}\,. \label{T2}
\end{eqnarray}
The estimates \eref{T1} and \eref{T2} imply that
\begin{equation}
\mathbb{E}\left(\sup_{ \|\xi\|=1}|\mathcal{A}_{2}|^p\right)
\leq C_{p}g(\varepsilon)^{2p}+C_{p}g(\varepsilon)^{p}\varepsilon^{2p}. \label{A2}
\end{equation}
Then by \eref{Mt},\eref{123},\eref{A12},\eref{A13} and \eref{A2}, for $\delta<1$, we obtain
\begin{eqnarray*}
\mathbb{P}\{\text{det} M_t(\underline{x})\leq \delta\}&\leq&
C_{p} g(\varepsilon)^{-p}(\varepsilon^{\eta_1 p}+\varepsilon^{\eta_2 p}+\varepsilon^{\eta_3 p}
+g(\varepsilon)^{2p}+g(\varepsilon)^{p}\varepsilon^{2p})\\
&\leq&C_{p}\delta^{\lambda p}\,,
\end{eqnarray*}
where $\lambda=\min\{\frac{\eta_1-\eta}{n\eta}, \frac{\eta_2-\eta}{n\eta}, \frac{\eta_3-\eta}{n\eta}, \frac{1}{n}, \frac{2}{n\eta}\}$.
The proof is completed.
\end{proof}

\section{Strict positivity of the density}
In this section, we proceed to the study of the positivity of the density $p_{t,\underline{x}}(\cdot)$ of the law of $u(t,\underline{x})$,
where $t>0$, $\underline{x}=(x_1,\ldots, x_n)$ are distinct points of $\RR^d$. 
The main theorem of this section is: 

\begin{theorem} \label{THM positive density}
Assume that conditions ({\bf H1})-({\bf H4}) hold, and the coefficients $\sigma$, $b$ are $\mathcal{C}^{\infty}$ functions with bounded derivatives of
all orders and $\sigma$ is bounded. We also assume $\sigma \neq 0$ on $\RR^d$.
Then the law of the random vector $u(t,\underline{x})$ has a $\mathcal{C}^{\infty}$ density $p_{t,\underline{x}}(y)$, and $p_{t,\underline{x}}(y)>0$ if $y$  belongs to the interior of the support
of the law of $u(t,\underline{x})$.
\end{theorem}

To prove this theorem we   will use the
criterion given  by Theorem 3.3 in \cite{BP}.  To state this
criterion in the context of framework, first we introduce  some notation and concepts.

Given    predictable processes $(g^1, \dots, g^n) \in \mathcal{H}^n_0$ and $z 
=(z_1, \cdots, z_n)\in \RR^n$, for any $h \in \mathcal{H}$ and $t\ge 0$
we define a  translation of $W_t(h)$:  
\begin{eqnarray*}
\widehat{W}_t (h):=\widehat{W}({\bf 1}_{[0,t]}h) = W ({\bf 1}_{[0,t]}h)+\sum_{k=1}^n z_k \langle {\bf 1}_{[0,t] }h, g^k\rangle_{\mathcal{H}_0}.
\end{eqnarray*}
Then $\{\widehat{W}_t, t\ge 0\} $ is a cylindrical Wiener process  in $\mathcal{H}$  on the probability space
$(\Omega, \mathcal{F}, \widehat{\mathbb{P}})$,
where
\begin{eqnarray*}
\frac{d \widehat {\mathbb{P}}}{d \mathbb{P}}  = \exp \left ( - \sum_{k=1}^n z_k \int_0^\infty \int_{\RR^d} g^k (s,y)W(ds,dy)-\frac{1}{2} \sum_{k=1}^n z_k^2 \int_0^\infty \| g^k (s,*)\|^2_{\mathcal{H}}ds\right)
 \,.
\end{eqnarray*}
Then, for any predicable process $Z \in L^2 (\Omega \times [0,\infty); \mathcal{H})$, we can write
\begin{eqnarray*}
&& \int_0^\infty \int_{\RR^d} Z(s,y)\widehat{W}(ds,dy)=\int_0^\infty\int_{\RR^d} Z(s,y)W(ds,dy)+ \sum_{k=1}^n z_k \int_0^\infty \langle Z(s,*), g^k(s,*)\rangle_{\mathcal{H}}
ds\,.
\end{eqnarray*}
For any $(t,x)\in [0,\infty)\times \RR^d$, let $\widehat {u}^z(t,x)$ be the solution to equation \eref{Eq} with respect to the cylindrical Wiener process $\widehat{W}$, that
is,
\begin{eqnarray}
\widehat{u}^z(t,x)&=&\int_0^t \int_{\RR^d} \Gamma(t-s, x-y)\sigma\left(\widehat{u}^z (s,y) \right)W(ds,dy) \nonumber\\
&&+\sum_{k=1}^n z_k \int_0^t \langle \Gamma(t-s, x-*)\sigma\left( \widehat {u}^z (s,*)\right), g^k(s,*)\rangle _{\mathcal{H}}ds \nonumber\\
&&+ \int_0^t \int_{\RR^d} b \left ( \widehat {u}^z(t-s,x-y)\right) \Gamma (s,dy)ds \label{u^z}\,.
\end{eqnarray}
Then, the law of $u$ under $\mathbb{P}$ coincides with the law of $\widehat{u}^z$ under $\widehat{\mathbb{P}}$.

Now we consider  a sequence $\{g_m\}_{m \geq 1}$ of predictable processes in $\mathcal{H}_0^n$ and $z \in \RR^n$. Let $\widehat{u}^z_m(t,x)$ be the solution to equation \eref{Eq}
with respect to the cylindrical Wiener process $\{\widehat{W}_t^m, t \ge 0\}$, where $\widehat{W}_t^m(h)=\widehat{W}^m({\bf 1}_{[0,t]}h)$ for any $h \in \mathcal{H}$,
and
\begin{eqnarray*}
\widehat{W}^m({\bf 1}_{[0,t]} h) = W ({\bf 1}_{[0,t]}h)+\sum_{k=1}^n z_k \langle {\bf 1}_{[0,t] }h, g_m^k\rangle_{\mathcal{H}_0} \,.
\end{eqnarray*}
Set $\varphi^z_{m,j}(t,x):= \partial_{z_j} \widehat{u}^z_m(t,x)$ and denote by $\varphi_m^z(t,\underline{x})$
 the $n \times n$ matrix $\{\varphi_{m,j}^z(t,x_i)\}_{1\leq i,j\leq n}$. Also, denote the Hessian matrix of $\widehat{u}^z_m(t,x)$ by $\psi_m^z(t,x):=
\partial_z^2 \widehat{u}_m^z(t,x)$, and let \newline $\psi_m^z(t,\underline{x}):=(\psi_m^z(t,x_1),\dots,\psi_m^z(t,x_n))$. In fact, it can be shown that
\begin{equation}
\partial_{z_j}\widehat{u}^z_m(t,x) = \int_0^t \langle D_{r,*} \widehat{u}^z_{m}(t,x), g_m^j(r,*) \rangle_{\mathcal{H}}dr\,.\nonumber
\end{equation}
We denote  the  operator norms of these matrices by $\|\varphi^z_m(t,x)\|$ and $\|\psi_m^z(t,x)\|$, respectively.

\medskip
 
We say that $y\in\RR^d$ satisfies ${\bf H}_{t, \underline{x}}(y)$ if there exists a sequence of predictable processes $\{g_{m}\}_{m\geq 1}$ in $\mathcal{H}_0^n$, and
positive constants
$c_1, c_2, r_0$ and $\delta$ such that

(i) $\limsup _{m \to \infty} \mathbb{P} \left\{ (\| u(t, \underline {x}) -y\| \leq r )
\cap ( | \text {det} \varphi_m^0(t, \underline{x}) | \geq c_1 )   \right\} > 0, \forall r \in (0, r_0]$.

(ii)
$\lim_{m \to \infty} \mathbb{P}\left \{ \sup_{|z|\leq \delta} \left( \|\varphi_m^z(t, \underline{x})\|
 + \| \psi_m^z (t, \underline{x})\| \right) \leq c_2  \right \}=1$.

Now we can state the criterion in \cite{BP} (Theorem 3.3) that we are going to use: 
Suppose that $y\in\RR^d$ belongs to the interior of the support of the law 
of    $u (t,\underline{x})$.  
If $y $ satisfies ${\bf H}_{t, \underline{x}}(y)$,  then $p_{t,\underline{x}}(y)>0$.  

\medskip 
{\noindent \it Proof  of Theorem \ref{THM positive density}}  \quad 
From the above criterion it suffices to check that $y$ satisfies the two conditions  in ${\bf H}_{t,\underline{x}}(y)$.  We will do this in several steps.

{\it Step 1}. Consider the sequence of predictable processes $\{g_m\}_{m\geq 1}$ in $\mathcal{H}^n_0$, defined by
\begin{eqnarray*}
g_m^k(s,*)=v_m^{-1}{\bf 1}_{[t-2^m, t]}(s)\Gamma(t-s, x_k-*)   \ \ \text{for $1 \leq k \leq n$}\,,
\end{eqnarray*}
where
\begin{eqnarray*}
v_m= \int_0^{2^{-m}}\int_{\RR^d} | \mathcal{F}\Gamma(r)(\xi)|^2 \mu(d\xi)dr\,.
\end{eqnarray*}
Taking the partial derivatives on both sides of  \eref{u^z} with $g$ replaced by $g_m$, we obtain that
\begin{eqnarray}
\partial_{z_j}\widehat{u}^z_m(t,x)&=&\int_{t-2^{-m}}^t \langle \Gamma(t-s, x-*)\sigma\left(\widehat {u}^z_m (s,*)\right), g^j_m(s,*) \rangle_{\mathcal{H}}ds \nonumber\\
&&+\sum_{k=1}^m z_k \int_{t-2^{-m}}^t \langle \Gamma(t-s, x-*)\sigma^{\prime}\left( \widehat{u}^z_m(s,*)\right)\partial_{z_j}\widehat{u}^z_m(s,*),
g^k_m(s,*)\rangle_{\mathcal{H}}ds \nonumber\\
&&+\int_{t-2^{-m}}^t \int_{\RR^d} \Gamma(t-s, x-y)\sigma^{\prime}\left( \widehat{u}^z_m(s,y)\right)\partial_{z_j}\widehat{u}^z_m(s,y)W(ds,dy) \nonumber\\
&&+ \int_{t-2^{-m}}^t \int_{\RR^d} b^{\prime} \left ( \widehat {u}^z_m (t-s,x-y)\right)\partial _{z_j}\widehat {u}^z_m(t-s,x-y)\Gamma(s,dy)ds  \nonumber \\
&:=&A_{m,j}^z(t,x)+B_{m,j}^z(t,x)+C_{m,j}^z(t,x)+D_{m,j}^z(t,x)\,.  \label{partial z}
\end{eqnarray}

{\it Step 2}. We are going to  bound  the   moments of the four terms on the right hand side of \eref{partial z}. We   assume that $\|z\|\leq \delta$ for some $\delta>0$.
Since $\sigma$ is bounded, there is a positive constant $K$ such that
\begin{equation}\label{E A^p}
|A_{m,j}^z(t,x)|\leq K \,.
\end{equation}
Using Minkowski's inequality and the fact that the partial derivatives of $\sigma$ are bounded, we get that for all $p\geq 1$, $t\leq T$,
\begin{equation}\label{E B^p}
\EE \left| B_{m,j}^z(t,x) \right|^p\leq C \delta^p \sup_{(s,y)\in [0,T]\times \RR^d}\EE \left | \partial _{z_j} \widehat {u}^z_m(s,y)\right|^p\,.\\
\end{equation} 
From the Burkholder-Davis-Gundy inequality and from the definition of $v_m$, we have 
\begin{eqnarray}
\EE \left |C_{m,j}^z(t,x)\right|^p
&\leq&C \sup_{(s,y)\in[0,T]\times \RR^d}\EE \left| \partial _{z_j} \widehat {u}^z_m (s,y)\right|^p  \left( \int_{t-2^{-m}}^t \int_{\RR^d} \left|
\mathcal{F}\Gamma(t-s)(\xi)\right|^2\mu(d\xi) ds\right)^{\frac{p}{2}} \nonumber\\
&\leq&C v_m^{\frac{p}{2}} \sup_{(s,y)\in[0,T]\times \RR^d} \EE \left| \partial_{z_j}\widehat{u}_m^z(s,y)\right|^p\,. \label{E C^p}
\end{eqnarray}
Since $b^{\prime}$ is bounded and by  condition ({\bf H1}),
\begin{equation}\label{E D^p}
\EE \left| D_{m,j}^z (t,x)\right|^p \leq C 2^{-mp} \sup_{(s,y)\in [0,T]\times \RR^d} \EE \left| \partial_{z_j}\widehat{u}_m^z(s,y)\right|^p\,.
\end{equation}
Combing \eref{E A^p}, \eref{E B^p}, \eref{E C^p} and  \eref{E D^p} we obtain
\begin{equation}\label{eq 4.1}
\sup_{(t,x)\in [0,T]\times \RR^d} \EE \left | \partial_{z_j}\widehat{u}_m^z (t,x)\right|^p\leq K + C (\delta^p+v_m^{\frac{p}{2}} + 2^{-mp}) \sup_{(s,y)\in
[0,T]\times \RR^d} \EE \left| \partial_{z_j}\widehat{u}_m^z(s,y)\right|^p\,.
\end{equation}
Proceeding as in the proof of Proposition 6.1 in \cite{NQ}, we can show 
\begin{eqnarray}\label{E partial u^z^p}
\sup_{(t,x)\in [0,T]\times \RR^d, |z|\leq \delta} \EE \left| \partial_{z_j}\widehat{u}^z_m (t,x)\right|^p < \infty\,.
\end{eqnarray}
Thus, when $m$ large enough, $\delta$ small enough, $C(\delta^p+v_m^{\frac{p}{2}}+2^{-mp})$ on the right hand side of equation \eref{eq 4.1} is less than $\frac{1}{2}$, we obtain
\begin{eqnarray}
\sup_{(t,x)\in [0,T]\times \RR^d,  |z|\leq \delta} \EE \left|\partial_{z_j} \widehat{u}_m^z(t,x) \right|^p \leq C\label{E partial u^z^p}
\end{eqnarray}
for some constant $C$.

Recall that $\varphi_{m,j}^z(t,x_i)=\partial_{z_j}\widehat{u}^z_m(t,x_i)$. Take $z=0$ and decompose $\varphi_{m,j}^0(t,x_i)$ as follows
\begin{equation}\label{varphi=ACD}
\varphi_{m,j}^0 (t,x_i)=A_{m,j}^0(t,x_i)+C_{m,j}^0(t,x_i) + D_{m,j}^0(t,x_i)\,.
\end{equation}
From  \eref{E C^p} and \eref{E D^p} it follows that
\begin{equation}\label{E CD^p}
\EE \left| C_{m,j}^0(t,x_i) +D_{m,j}^0(t,x_i) \right|^p \leq C (v_m^{\frac{p}{2}}+2^{-mp})\,.
\end{equation}
For $A_{m,j}^0(t,x_i)$,
\begin{eqnarray}
A_{m,j}^0(t,x_i)&=&\int_{t-2^{-m}}^t \langle \Gamma(t-s, x_i-*)\sigma\left( u(s,*)\right), g_m^j(s,*)\rangle_{\mathcal{H}}ds \nonumber\\
&=& \int_{t-2^{-m}}^t \langle \Gamma(t-s,x_i-*) \left [ \sigma \left( u(s,*) \right) - \sigma \left( u(t,x_i) \right)\right],
g_m^j(s,*)\rangle_{\mathcal{H}}ds \nonumber\\
&&+ \sigma\left (u(t,x_i) \right) \int_{t-2^{-m}}^t \langle \Gamma (t-s,x_i-*), g_m^j(s,*)\rangle_{\mathcal{H}}ds \nonumber\\
&:=& O_{m,i,j}+ \widetilde{O}_{m,i,j}\,. \label{A=O+tilde O}
\end{eqnarray}
By the assumption ({\bf H2}) and Minkowski's inequality, we have
\begin{eqnarray*}
&&\EE \left| O_{m,i,j}\right|^p\\
 &=& \left \| \frac{1}{v_m} \int_{t-2^{-m}}^t  \int_{\RR^d} \int_{\RR^d}\Gamma(t-s, x_i-dy) \left[ \sigma\left(
u(s,y)\right)-\sigma\left( u(t,x_i)\right)\right] \Gamma(t-s,x_j-dz)f(y-z)ds\right\|^p_{L^p(\Omega)}\\
&\leq&\frac{1}{v_m^p} \left( \int_{t-2^{-m}}^t \int_{\RR^d} \int_{\RR^d} \| \sigma\left( u(s,y)\right)-\sigma\left( u(t,x_i)\right)\|_{L^p(\Omega)}
\Gamma(t-s,x_j-dz)f(y-z)\Gamma(t-s,x_i-dy) ds\right)^p\\
&\leq& \frac{C}{v_m^p} \left ( \int_{t-2^{-m}} ^ t \int_{\RR^d} \int_{\RR^d} \Gamma(t-s,x_i-dy) \left( |x_i-y|^{\kappa_2} + |s-t|^{\kappa_1}\right)
\Gamma(t-s,x_j-dz)f(y-z)ds\right)^p\\
&\leq&\frac{C}{v_m^p} \left( 2^{-m\eta_1}+2^{-m\eta_3}\right)^p \to 0 \ \text{as $ m \to \infty$. }
\end{eqnarray*}
For $\widetilde{O}_{m,i,j}$, when $i=j$, it is easy to see that
\begin{equation}\label{tilde O}
 \widetilde{O}_{m,i,i}=\sigma(u(t,x_i))\,,
\end{equation}
 while when $i \neq j$, we have the $p$th moment bound
\begin{eqnarray}
\EE \left| \widetilde{O}_{m,i,j}\right|^p&\leq& \EE \left| \sigma (u(t,x_i))\right|^p \left( \int_{t-2^{-m}}^t \langle \Gamma(t-s,x_i-*), g_m^j(s,*)\rangle_{\mathcal{H}}ds\right)^p \nonumber\\
&\leq&C_p \left(\frac{1}{v_m}\int_0^{2^{-m}}\int_{\RR^d}\int_{\RR^d}\Gamma(s,x_i-dy)f(y-z)\Gamma(s,x_j-dz)ds\right)^p \nonumber\\
&\leq &C_p \left( \frac{2^{-m\eta_2}}{v_m}\right)^p\,, \label{EtildeO^p}
\end{eqnarray}
which goes to $0$ as $m \to \infty$.

{\it Step 3}.  We check condition (i) in  hypothesis ${\bf H}_{t,\underline{x}}(y)$.
Recall that $y \in \text{Supp} \left(P_{u(t,\underline{x})}\right)$, there exists $r_0>0$ such that for all $0 < r \leq r_0$,
\begin{equation}
\mathbb{P}\left\{u(t,\underline{x}) \in B(y;r) \right\} > 0\,.\nonumber
\end{equation}
By the assumption on $\sigma$, there is a $c_1 > 0$ such that
\begin{equation}\label{P sigma}
\mathbb{P}\left \{  \left (\| u(t,\underline{x})-y\| \leq r \right) \cap \left ( \left| \prod_{i=1}^n \sigma (u(t,x_i))\right| \geq 2c_1\right) \right \} > 0
\end{equation}
where
\begin{equation}
c_1 =\frac{1}{2} \left(\inf_{z \in B(y; r) } |\sigma(z)|\right)^n\,.\nonumber
\end{equation}
Recall that  $\varphi_m^0(t,\underline{x})$ is the matrix $\left ( \varphi_{m,j}^0 (t,x_i)\right) _{1 \leq i, j\leq n}$, by \eref{varphi=ACD}, \eref{E CD^p}, \eref{A=O+tilde O}, \eref{tilde O} and  \eref{EtildeO^p}, we obtain

\begin{equation}\label{E det error}
\EE \left| \text{det} \varphi_m^0(t,\underline{x}) - \prod_{i=1}^n \sigma (u(t,x_i))\right|^p \to 0 \ \text{as}\  m \to \infty\,.
\end{equation}
Combination of  \eref{P sigma} and \eref{E det error}  yields
\begin{eqnarray*}
\limsup_{m \to \infty} \mathbb{P} \left \{ \left( \| u(t, \underline{x}) -y\| \leq r \right) \cap \left( \left |\text{det} \varphi_m^0(t,\underline{x}) \right| \geq c_1\right) \right \} > 0 \,.
\end{eqnarray*}

{\it Step 4}. We check condition (ii) in the hypothesis ${\bf H}_{t,\underline{x}}(y)$.

We first show that there exists $c_2> 0$ and $\delta > 0$ such that
\begin{eqnarray*}
\lim_{m \to \infty} \mathbb{P} \left\{ \sup_{|z|\leq \delta} \| \varphi_m^z(t,\underline{x})\| \leq c_2  \right\} =1 \,.
\end{eqnarray*}
Consider the following equation
\begin{equation}\label{eq v^z}
v_{m,j}^z(t,x)=A_{m,j}^z(t,x)+\sum_{k=1}^n z_k \int_{t-2^{-m}}^t \langle \Gamma(t-s, x-*) \sigma^{\prime}\left(\widehat {u}_m^z(s,*)\right)v_{m,j}^z(s,*), g^k(s,*)
\rangle_{\mathcal{H}}ds \,. 
\end{equation}
By the  contraction mapping theorem we can prove that this equation has a unique solution $v_{m,j}^z(t,x)$ and there exists a constant $C$ such that
\begin{equation}\label{v^z bdd}
\sup_{(t,x)\in [0,T]\times \RR^d, |z|\leq
\delta} | v_{m,j}^z (t,x)|\leq C \quad \forall 1 \leq j \leq n\,,
\end{equation}
when $\delta$ is small.

Then we claim that for each $j$, $v_{m,j}^z(t,x)-\partial_{z_j}\widehat{u}_m^z(t,x)$ converges to $0$ in $L^p(\Omega)$ norm, uniformly in $(t,x)\in [0,T]\times \RR^d$, and $|z|\leq \delta$ when $\delta$ is small.  Indeed,  we have 
\begin{eqnarray*}
&&\EE \left| \partial_{z_j} \widehat{u}_m^z(t,x)-v_{m,j}^z(t,x)\right|^p\\
&\leq&C_p \sum_{k=1}^n |z_k|^p \left( \int_{t-2^{-m}}^t \langle \Gamma(t-s,x-*)\sigma^{\prime}(\widehat{u}_m^z(s,*)), g^k(s,*) \rangle_{\mathcal{H}} ds \right)^p \\
&&\quad \times\sup_{(s,y)\in [0,T]\times \RR^d} \EE \left| \partial_{z_j}\widehat{u}^z_m(s,y)-v_{m,j}^z(s,y) \right|^p\\
&&+C_p \left\| \int_{t-2^{-m}}^t \int_{\RR^d} \Gamma(t-s,x-y)\sigma^{\prime}\left( \widehat{u}_m^z(s,y)\right)\partial_{z_j}\widehat{u}_m^z(s,y)W(ds,dy)
\right\|_{L^p(\Omega)}^p\\
&&+C_p \left\| \int_{t-2^{-m}}^t \int_{\RR^d} b^{\prime} \left( \widehat{u}_m^z (t-s, x-y)\right) \partial_{z_j} \widehat{u}_m^z(t-s, x-y) \Gamma(s,dy)ds\right
\|_{L^p(\Omega)}^p\\
&\leq& C_p \delta^p \sup_{(t,x)\in [0,T]\times \RR^d} \EE \left | \partial_{z_j} \widehat{u}_m^z(t,x)-v_{m,j}^z(t,x) \right |^p\\
&&+ C_p \left ( \int_{t-2^{-m}}^t \int_{\RR^d} \int_{\RR^d} \Gamma(t-s,x-dy)f(y-\tilde{y}) \Gamma(t-s, x-d\tilde{y}) ds \right)^{\frac{p}{2}} \\
&&\quad \times \sup_{(s,y)\in[0,T] \times \RR^d} \EE \left | \partial_{z_j} \widehat{u}_m^z(s,y) \right|^p\\
&&+ C_p \left( \int_{t-2^{-m}}^t \int_{\RR^d} \Gamma(s,dy)ds \right)^p \sup_{(s,y)\in [0,T]\times \RR^d}\EE \left |\partial _{z_j} \widehat{u}_m ^z
(s,y)\right|^p\\
&\leq& C_p \delta^p \sup_{(t,x)\in [0,T]\times \RR^d} \EE \left | \partial_{z_j} \widehat{u}_m^z(t,x)-v_{m,j}^z(t,x) \right |^p\\
&&+ C_p \left ( \int_0^{2^{-m}} \int_{\RR^d} |\Gamma(s)(\xi)|^2 \mu(d\xi) ds \right)^{\frac{p}{2}}
 \sup_{(s,y)\in[0,T] \times \RR^d} \EE \left | \partial_{z_j} \widehat{u}_m^z(s,y) \right|^p\\
&&+ C_p \left( \int_{t-2^{-m}}^t \int_{\RR^d} \Gamma(s,dy)ds \right)^p \sup_{(s,y)\in [0,T]\times \RR^d}\EE \left |\partial _{z_j} \widehat{u}_m ^z
(s,y)\right|^p\,. 
\end{eqnarray*}
Thus  we can choose $\delta$ small enough such that 
$C_p\delta^p \leq \frac{1}{2}$, using condition ({\bf H1}) 
and \eref{E partial u^z^p} to conclude that
\begin{equation}\label{u^z -v^z}
\sup_{(t,x)\in [0,T]\times \RR^d, |z|\leq \delta} \EE \left | \partial_{z_j} \widehat{u}^z_m(t,x)- v_{m,j}^z (t,x)\right|^p
\end{equation}
goes to $0$ as $m$ tends to $\infty$.

Next, we will calculate the $p$th moment of the increments with respect to $z$ of $\partial_{z_j}\widehat{u}_m^z(t,x)$ and $v_{m,j}^z(t,x)$.
\begin{eqnarray*}
&&\EE \left| \partial_{z_j}\widehat{u}_m^z (t,x) - \partial_{z_j} \widehat{u}^{z^{\prime}}_m (t,x)  \right|^p\\
&\leq& \EE \left |\int_{t-2^{-m}}^t \left\langle  \Gamma(t-s, x-*) \left [ \sigma \left( \widehat{u}_m^z(s,*) \right)- \sigma ( \widehat{u}_m^{z^{\prime}}(s,*)
)\right], g^j_m(s,*)\right\rangle _{\mathcal{H}} ds \right | ^p\\
&&+\EE \Big | \sum_{k=1}^n \int_{t-2^{-m}}^t \langle \Gamma(t-s,x-*) [z_k \sigma^{\prime} \left (\widehat{u}_m^z(s,*) \right) \partial_{z_j}\widehat{u}_m^z(s,*)\\
 && \quad -z^{\prime}_k \sigma^{\prime} (\widehat{u}_m^{z^{\prime}}(s,*) ) \partial_{z_j}\widehat{u}_m^{z^{\prime}}(s,*)], g^k_m(s,*)\rangle_{\mathcal{H}} ds
\Big|^p\\
&&+ \EE \left| \int _{t-2^{-m}}^t \int_{\RR^d} \Gamma(t-s,x-y) [\sigma^{\prime} \left (\widehat{u}_m^z(s,y) \right) \partial_{z_j}\widehat{u}_m^z(s,y)-
\sigma^{\prime} (\widehat{u}_m^{z^{\prime}}(s,y) ) \partial_{z_j}\widehat{u}_m^{z^{\prime}}(s,y)] W(ds,dy)\right |^p \\
&&+ \EE \Big | \int_{t-2^{-m}} ^t \int_{\RR^d} \Big[b^{\prime}\left ( \widehat{u}_m^z (t-s,x-y)\right)\partial_{z_j}\widehat{u}_m^z(t-s,x-y)\\
&& \quad - b^{\prime}( \widehat{u}_m^{z^{\prime}} (t-s,x-y))\partial_{z_j} \widehat{u}_m^{z^{\prime}}(t-s,x-y)\Big] \Gamma(s,dy) ds \Big|^p\,.
\end{eqnarray*}
Proceeding as before, we obtain that
\begin{eqnarray*}
\EE \left| \partial_{z_j}\widehat{u}_m^z (t,x) - \partial_{z_j} \widehat{u}^{z^{\prime}}_m (t,x)  \right|^p \leq C |z-z^{\prime}|^p
\end{eqnarray*}
uniformly in $(t,x)\in [0,T]\times \RR^d$, $|z|\leq \delta$ and $m$.  Similarly, we have
\begin{eqnarray*}
\EE \left|v_{m,j}^z(t,x)- v_{m,j}^{z^{\prime}}(t,x)\right|^p\leq C |z-z^{\prime}|^p
\end{eqnarray*}
uniformly in $(t,x)\in [0,T]\times \RR^d$, $|z|\leq \delta$ and $m$.  Using Kolmogorov's continuity theorem and \eref{v^z bdd}, \eref{u^z -v^z} we obtain
\begin{equation}
\lim_{m \to \infty} \mathbb{P}\left \{ \sup_{|z|\leq \delta} \|\varphi_m^z(t, \underline{x})\| \leq C  \right \}=1\nonumber
\end{equation}
for some positive constant $C$.

Next we will show that there exists a positive constant $C$ such that
\begin{equation*}
\lim_{m \to \infty}\mathbb{P}\left\{\sup_{|z|\leq \delta} \| \psi_m^z(t,\underline{x})\|\leq C \right\}=1\,.
\end{equation*}
This proof is analogous to that for $\varphi_m^z(t,\underline{x})$, but the computations are more involved. Let us just write the equation for the quantity of interest and the main steps.  Take the partial derivative on  both sides of  \eref{partial z},  we obtain
\begin{eqnarray*}
&&\partial_{z_l}\partial_{z_j}\widehat{u}^z_m(t,x)\\
&=&\int_{t-2^{-m}}^t \langle \Gamma(t-s, x-*)\sigma^{\prime}\left(\widehat {u}^z_m (s,*)\right)\partial_{z_l}\widehat{u}_m^z(s,*), g^j_m(s,*) \rangle_{\mathcal{H}}ds\\
&&+ \int_{t-2^{-m}}^t \langle \Gamma(t-s, x-*)\sigma^{\prime}\left( \widehat{u}^z_m(s,*)\right)\partial_{z_j}\widehat{u}^z_m(s,*),
g^l_m(s,*)\rangle_{\mathcal{H}}ds\\
&&+\sum_{k=1}^m z_k \int_{t-2^{-m}}^t \langle \Gamma(t-s, x-*)\Big(\sigma^{\prime\prime}\left( \widehat{u}^z_m(s,*)\right)\partial_{z_l}\widehat{u}^z_m(s,*)\partial_{z_j}\widehat{u}^z_m(s,*)\\
&&\quad +\sigma^{\prime}\left( \widehat{u}^z_m(s,*)\right)\partial_{z_l}\partial_{z_j}\widehat{u}^z_m(s,*)\Big), g^k_m(s,*)\rangle_{\mathcal{H}}ds\\
&&+\int_{t-2^{-m}}^t \int_{\RR^d} \Gamma(t-s, x-y)        \Big(\sigma^{\prime\prime}\left( \widehat{u}^z_m(s,y)\right)\partial_{z_l}\widehat{u}^z_m(s,y)\partial_{z_j}\widehat{u}^z_m(s,y)\\
&&\quad +\sigma^{\prime}\left( \widehat{u}^z_m(s,y)\right)\partial_{z_l}\partial_{z_j}\widehat{u}^z_m(s,y)\Big)     W(ds,dy)\\
&&+ \int_{t-2^{-m}}^t \int_{\RR^d} \Big(b^{\prime\prime}\left( \widehat{u}^z_m(t-s,x-y)\right)\partial_{z_l}\widehat{u}^z_m(t-s,x-y)\partial_{z_j}\widehat{u}^z_m(t-s,x-y)\\
&&\quad +b^{\prime}\left( \widehat{u}^z_m(t-s,x-y)\right)\partial_{z_l}\partial_{z_j}\widehat{u}^z_m(t-s,x-y)\Big) \Gamma(s,dy)ds
\end{eqnarray*}
 and a similar equation for $\partial_{z_l}v_{m,j}^z(t,x)$.

 We can show that for every $1 \leq l,j \leq n$,
 \begin{equation*}
\sup_{(t,x)\in [0,T]\times \RR^d, |z|\leq \delta } \EE \left|\partial_{z_l}\partial_{z_j}\widehat{u}_m^z(t,x)-\partial_{z_l}v_{m,j}^z(t,x)\right|^p \to 0\,,
 \end{equation*}
 as $m$ goes to $\infty$. Bound $\partial_{z_l}v_{m,j}^z(t,x)$ and calculate the $p$th moment of the increments with respect to $z$ of $\partial_{z_l}\partial_{z_j}\widehat{u}_m^z(t,x)$ and $\partial_{z_l}v_{m,j}^z(t,x)$. The result follows as in the previous step.

{\it Step 5}.  Combing the results in step 3 and step 4, together with the criterion developed by Theorem 3.3 in \cite{BP} that we cited just before the proof, we complete the proof. 
 \rule{0.5em}{0.5em}

\section{Examples}
In this section  we will give some examples of  fundamental solutions $\Gamma$ and  covariance functions $f$ satisfying hypotheses  ${(\bf H1)}$ to ${(\bf H4)}$. This implies that Theorem \ref{THM smooth density} and Theorem \ref{THM positive density} can be applied to these examples. We consider the fundamental solution to the heat  equation in any dimension and the  wave equation in dimensions up to three and the covariance functions given by the Riesz, Bessel, and fractional kernels.

\subsection{Heat equation}
Let $\Gamma(r,dx)$ be the fundamental solution to the heat equation on $\RR^d$, i.e., $\Gamma(r,dx)=p_r(x)dx$, where $p_r(x)=(2\pi r)^{-d/2}e^{-\frac{|x|^2}{2r}}$ is the $d$-dimensional heat kernel. Then, hypothesis ${(\bf H1)}$ to ${(\bf H4)}$ are satisfied for the following covariance functions:

\medskip
\noindent
(A) \textit {Riesz kernel}. Let  $f(x)=|x|^{-\beta}$ with $0 < \beta < 2 \wedge d$.
It is well-known that {(\bf H1)} holds.  According to \cite{SS},  $({\bf H2})$ is satisfied with $0< \kappa_1<\frac{2-\beta}{4}$ and $0< \kappa_2 < \frac{2-\beta}{2}$.   In \cite{NuaE} it is proved that
 $({\bf H3})$  holds with $\eta = \frac{2-\beta}{2}$, and property (i) in (${\bf H 4}$) holds with  $\eta_1= \frac {2-\beta} 2 +\kappa_1$.

  Next we check conditions (ii) and (iii) in (${\bf H 4}$).    To show  \eref{hyp 4-1} we use the fact that  there exists a constant $C>0$ such that for  any non zero $y\in \mathbb{R}^d$ and $r\ge 0$
  \begin{equation} \label{5.1}
  \int_{\mathbb{R}^d} p_r(x) |x-y|^{-\beta} dx \le C |y|^{-\beta}.
  \end{equation}
For a non zero $w\in \RR^d$,  using (\ref{5.1}) we can write
\begin{eqnarray*}
 \int_0^{\varepsilon}\langle p_r (*), p_r(w+*)\rangle_{\mathcal{H}}dr
&=&\int_0^{\varepsilon} \int_{\RR^d} \int_{\RR^d} p_r(x)p_r(y+w)|x-y|^{-\beta}dx dy dr \\
&\leq&C \int_0^{\varepsilon} \int_{\RR^d} p_{ r} (y+w) |y|^{-\beta} dy dr \\
& \leq& C   \varepsilon  |w|^{-\beta},
\end{eqnarray*}
so \eref{hyp 4-1} is satisfied with $\eta_2=1 > \eta$.  For \eref{hyp 4-3}, using the fact that $\sup_{x \in \RR^d}|x|^{\alpha}e^{-x^2}< \infty$ for any positive $\alpha$, we have
\begin{eqnarray*}
 \int_0^{\varepsilon}\langle |*|^{\kappa_2}p_r(*), p_r(w+*) \rangle_{\mathcal{H}}dr
&=&\int_0^{\varepsilon} \int_{\RR^d}\int_{\RR^d} |x|^{\kappa_2}p_r(x)p_r(y+w)|x-y|^{-\beta} dx dy dr \\
&\leq&C \int_0^{\varepsilon} \int_{\RR^d} \int_{\RR^d} r^{\frac{\kappa_2}{2}} p_{2r}(x)p_r(y+w)|x-y|^{-\beta}dx dy dr\\
&\leq& C \int_0^{\varepsilon}r^{\frac{\kappa_2}{2}} \int_{\RR^d} e^{-r |\xi|^2} e^{-\frac{1}{2}r |\xi|^2} |\xi|^{\beta-d}d\xi dr \\
&=& C \int_0^{\varepsilon} r^{\frac{\kappa_2-\beta}{2}}dr =  C \varepsilon^{\frac{\kappa_2-\beta}{2}+1}\,.
\end{eqnarray*}
Therefore,  \eref{hyp 4-3} is satisfied with $\eta_3=\frac{\kappa_2-\beta}{2}+1> \eta$.

\medskip
\noindent
(B)  \textit{Bessel  kernel}.
Let $f(x)=\int^{\infty}_{0}u^{\frac{\alpha-d-2}{2}}e^{-u}e^{-\frac{|x|^2}{4u}}du$, $ d-2<\alpha<d$.   In this case $\mu(d\xi)= c_{\alpha,d} (1+ |\xi|^2) ^{-\frac {\alpha}2 } d\xi$.
Hypothesis {(\bf H1)} can be easily verified by direct computation.  According to \cite{SS},  $({\bf H2})$ is satisfied with $0< \kappa_1<\frac{2-d+\alpha}{4}$ and $0< \kappa_2 < \frac{2-d+\alpha}{2}$.  For $({\bf H3})$, we note that, assuming $\varepsilon <1$,
\begin{eqnarray*}
\int_0^{\varepsilon} \int_{\RR^d} |\mathcal{F}\Gamma(r)(\xi)|^2\mu(d\xi)dr&=&C \int_0^{\varepsilon} \int_{\RR^d}e^{-r|\xi|^2}  (1+ |\xi|^2)  ^{-\frac   {\alpha}{2}} d\xi dr \\
&=&C \int_0^{\varepsilon}\int_{\RR^d}e^{-|\theta|^2}\frac{r^{\frac{\alpha-d}{2}}}{(|\theta|^2+ r)^{ \frac {\alpha}2 }}d\theta dr\\
&\geq& C \int_0^{\varepsilon}r^{\frac{\alpha-d}{2}}dr\int_{\RR^d}e^{-|\theta|^2}\frac{1}{(|\theta|^2+1)^{ \frac {\alpha}2}}d\theta\\
&=& C \varepsilon^{\frac{\alpha-d}{2}+1}\,.
\end{eqnarray*}
Thus,  $({\bf H3})$ is satisfied with $\eta = \frac{\alpha-d}{2}+1$. To show (${\bf H 4}$) we use the fact that  for any $x \in \RR^d$, $f(x)\leq C |x|^{-d+\alpha}$ (see Proposition 6.1.5 in \cite{Gra}). Therefore, proceeding as in the case of the Riesz kernel with $\beta= d-\alpha$ we obtain that
  conditions \eref{hyp 4-2}, \eref{hyp 4-1} and  \eref{hyp 4-3} in (${\bf H 4}$)  hold, with $\eta_1 = \frac {\alpha -d} 2 +1 +\kappa_1$,
    $\eta_2=1 $ and  $\eta_3 = \frac {\alpha -d} 2 +1 +\frac {\kappa_2}2$, respectively.

\medskip
\noindent
(C) {\textit{Fractional kernel}.
Let $f(x)=\prod_{j=1}^d |x_j|^{2H_j-2}$, $ \frac{1}{2} < H_j < 1$ for $1\leq j \leq d$ such that $\sum_{j=1}^d H_j > d-1$.  First notice that although we have assumed $f(x)$ to be a continuous function on $\RR^d \setminus \{0\}$, it is clear that all of our theory still works for this case. Then we note that since $f(x)=\prod_{j=1}^d |x_j|^{2H_j-2}$, we have $\mu(d\xi)=C_H\prod_{j=1}^d |\xi_j|^{1-2H_j}d\xi$, where $C_H$ only depends on $H:=(H_1,H_2,\dots,H_d)$.
According to \cite{SS},  ({\bf H1})  holds and ({\bf H2}) is satisfied for $0 < \kappa_1 < \frac{1}{2}(\sum_{j=1}^d H_j -d+1)$ and $0 < \kappa_2 < \sum_{j=1}^d H_j -d+1$.   For ({\bf H3}), using the change of variable $\sqrt{t}\xi \rightarrow \xi$, we obtain
 \begin{eqnarray*}
 \int_{0}^\varepsilon \int_{\RR^d} |\mathcal{F}\Gamma(t)(\xi)|^2\mu(d\xi)dt
 &=&\int_0^\varepsilon \int_{\RR^d} e^{-t|\xi|^2}\prod_{j=1}^d |\xi_j|^{1-2H_j}d\xi dt =C \varepsilon^{\sum_{j=1}^d H_j -d+1},
  \end{eqnarray*}
so ({\bf H3}) is verified with $\eta = \sum_{j=1}^d H_j-d+1$.
For \eref{hyp 4-2},  we can proceed as in checking ({\bf H3}) to get
\begin{eqnarray*}
\int_0^{\varepsilon} r^{\kappa_1} \|\Gamma(r)\|_{\mathcal{H}}^2 dr = C \int_0^{\varepsilon} r^{\kappa_1 + \sum_{j=1}^d H_j-d} dr = C \varepsilon ^{\sum_{j=1}^d H_j-d+1+\kappa_1}\,,
\end{eqnarray*}
so \eref{hyp 4-2} is satisfied with $\eta_2 = \sum_{j=1}^d H_j-d+1+\kappa_1$ which is strictly greater than $\eta$.

 To check \eref{hyp 4-1}, fix  a nonzero point $w = (w_1, w_2,\dots,w_d) \in \RR^d$, without loss of generality, we may assume that $w_1 \neq 0$. Then using Fourier transform and \eref{5.1} we have
\begin{eqnarray*}
&&\int_0^{\varepsilon} \langle \Gamma(r,*), \Gamma(r,w+*)\rangle_{\mathcal{H}}dr
 = \int_0^{\varepsilon}\int_{\RR^d}\int_{\RR^d} p_r(x)  p_r(w+y)\prod_{j=1}^d |x_j-y_j|^{2H_j-2} dy dx dr \\
&=&\int_0^{\varepsilon} \prod_{j=1}^d\left (\int_{\mathbb{R}}\frac{1}{(2\pi r)^{\frac{1}{2}}}e^{-\frac{|x_j|^2}{2r}}\frac{1}{(2\pi r)^{\frac{1}{2}}}e^{-\frac{|w_j+y_j|^2}{2r}} |x_j-y_j|^{2H_j-2} dy_j dx_j \right) dr \\
&=&C \int_0^{\varepsilon} \left (\int_{\mathbb{R}}\frac{1}{(2\pi r)^{\frac{1}{2}}}e^{-\frac{|x_1|^2}{2r}}\frac{1}{(2\pi r)^{\frac{1}{2}}}e^{-\frac{|w_1+y_1|^2}{2r}} |x_1-y_1|^{2H_1-2} dy_1 dx_1 \right) \\
 && \times \prod_{j=2}^d \left ( \int_{\mathbb{R}} e^{-r |\xi_j|^2} e^{-i w_j \xi_j} |\xi_j|^{1-2H_j} d\xi_j\right)dr \\
&\leq& C \int_0^{\varepsilon} \left (\int_{\mathbb{R}}\frac{1}{(2\pi r)^{\frac{1}{2}}}e^{-\frac{|x_1|^2}{2r}}\frac{1}{(2\pi r)^{\frac{1}{2}}}e^{-\frac{|w_1+y_1|^2}{2r}} |x_1-y_1|^{2H_1-2} dy_1 dx_1 \right) \\
 && \times \prod_{j=2}^d \left ( \int_{\mathbb{R}} e^{-r |\xi_j|^2}  |\xi_j|^{1-2H_j} d\xi_j\right)dr \\
&\leq&C |w_1|^{2H_1-2} \int_0^{\varepsilon} r^{\sum_{j=2}^d H_j-d+1} dr = C \varepsilon^{\sum_{j=2}^d H_j -d+2}\,,
\end{eqnarray*}
where in the last inequality we have used the change of variable $\sqrt{r}\xi \rightarrow \xi$.  So \eref{hyp 4-1} is satisfied with $\eta_1=\min_{1\leq k \leq d}(\sum_{j\neq k}^d H_j -d +2)$, which is strictly greater than $\eta$.   For \eref{hyp 4-3}, fixing again a non zero element $w\in \mathbb{R}^d$ and using the bound $|x|^{\alpha}p_r(x)\leq C r^{\frac{\alpha}{2}} p_{2r}(x)$,  for all $  x \in \RR^d$,   we have
\begin{eqnarray*}
  \int_0^{\varepsilon} \langle |*|^{\kappa_2} \Gamma(r,*), \Gamma(r,w+*) \rangle_{\mathcal{H}}dr
 &=& \int_0^{\varepsilon} \int_{\RR^d}\int_{\RR^d} |x|^{\kappa_2} p_r(x) p_r(y+w)  \prod_{j=1}^d |x_j-y_j|^{2H_j-2}dx dy dr\\
&\leq&C \int_0^{\varepsilon} \int_{\RR^d}\int_{\RR^d} r^{\frac{\kappa_2}{2}}  p_{2r}(x) p_r(y+w) \prod_{j=1}^d |x_j-y_j|^{2H_j-2}dx dy dr\\
&\leq&C \int_0^{\varepsilon} \int_{\RR^d} r^{\frac{\kappa_2}{2}} e^{-\frac{3r}{2} |\xi|^2} \prod_{j=1}^d |\xi_j|^{1-2H_j} d\xi dr\\
&=&C \int_0^{\varepsilon} r^{\frac{\kappa_2}{2}+\sum_{j=1}^d H_j -d}dr= C \varepsilon^{\frac{\kappa_2}{2}+\sum_{j=1}^d H_j -d+1} \,,
\end{eqnarray*}
so \eref{hyp 4-3} is satisfied with $\eta_3=\frac{\kappa_2}{2}+\sum_{j=1}^d H_j -d+1$, which is strictly greater than $\eta$.

\subsection{Wave equation}
Let $\Gamma_{d}(t,dx)$ be the fundamental solution to the wave equation on $\RR^d$, for $d=1,2,3$, i.e., $\Gamma_1(t,dx)=\frac{1}{2}{\bf 1}_{\{|x|<t\}}dx$, $\Gamma_2(t,dx)=\frac{1}{2\pi}(t^2-|x|^2)^{-1/2}_{+}dx$, $\Gamma_3(t,dx)=\frac{1}{4\pi t}\sigma_t(dx)$, where $\sigma_t$ denotes the surface measure on the two-dimensional sphere of radius $t$. We recall  that the Fourier transform of $\Gamma_d(t,dx)$ is given by
$$
\mathcal{F}\Gamma_d(t)(\xi)=\frac{\text{sin}( t|\xi|)}{ |\xi|}.
$$

\medskip
\noindent
(A) \textit {Riesz kernel}.
 Let  $f(x)=|x|^{-\beta}$ with $0 < \beta < 2\wedge d$. It is known that hypothesis ${(\bf H1)}$is satisfied.
   According to \cite{HHN},  $({\bf H2})$ is satisfied with $0< \kappa_1=\kappa_2 < \frac{2-\beta}{2}$.
 In \cite{NuaE} it is proved that condition {(\bf H3)} is satisfied for $\eta=3-\beta$ and   \eref{hyp 4-2} holds with
 $\eta_1= \kappa_1+3-\beta >\eta$.
  To show \eref{hyp 4-1}, we fix $w\neq0$, and taking $\varepsilon$   such that $4\varepsilon<|w|$ we get
$\frac{|w|}{2}\leq |x-y|\leq \frac{3|w|}{2}$ if $|x|\leq  \varepsilon$ and $|w+y| \leq \varepsilon$. Then, $|x-y|^{-\beta}$ is bounded by some constant $C$ depending on $|w|$. Hence we have
\begin{eqnarray*}
 \int_0^{\varepsilon} \langle \Gamma_d(r,*), \Gamma_d(r,w+*)\rangle_{\mathcal{H}}dr
&=&\int_0^{\varepsilon} \int_{\RR^d} \int_{\RR^d} \Gamma_d(r, dx)\Gamma_d(r, w+dy)|x-y|^{-\beta}dr \\
&\leq&C_w \int_0^{\varepsilon} \int_{\RR^d} \Gamma_d(r,dx) \int_{\RR^d}\Gamma_d(r,w+dy)dr\\
&\leq&C_{w} \int_0^{\varepsilon}r^2 dr \leq C_{w}\varepsilon^3 \,,
\end{eqnarray*}
so \eref{hyp 4-1} is satisfied with $\eta_2=3 > \eta$.
For \eref{hyp 4-3}, any fixed $w\in\RR^d$, using again the same arguments, we have
\begin{eqnarray*}
 \int_0^{\varepsilon} \langle  |*|^{\kappa_2} \Gamma_d(r,*), \Gamma_d(r,w+*)\rangle_{\mathcal{H}}dr
&=&\int_0^{\varepsilon} \int_{\RR^d} \int_{\RR^d} |x|^{\kappa_2}\Gamma_d(r,dx) \Gamma_d(r,w+dy)|x-y|^{-\beta}dr\\
&\leq& \int_0^{\varepsilon} |r|^{\kappa_2}\int_{\RR^d} \int_{\RR^d}\Gamma_d(r,dx) \Gamma_d(r,w+dy)|x-y|^{-\beta}dr\\
& \leq & C  \varepsilon ^{\kappa_2 +3-\beta},
\end{eqnarray*}
so \eref{hyp 4-3} is satisfied with $\eta_3=\kappa_2+3-\beta> \eta$.

\medskip
\noindent
(B) \textit {Bessel kernel}.
 Let $f(x)=\int^{\infty}_{0}u^{\frac{\alpha-d-2}{2}}e^{-u}e^{-\frac{|x|^2}{4u}}du$, $ \max(d-2,0)<\alpha<d$.   According to section 3 in \cite{NuaE} and  \cite{HHN}, {(\bf H1)} holds and   $({\bf H2})$ is satisfied with $0< \kappa_1=\kappa_2 < \frac{\alpha-d+2}{2}$.
Making the change of variable $r\xi \rightarrow \xi$ and  assuming $\varepsilon<1$, we get that
\begin{eqnarray*}
\int_0^{\varepsilon} \int_{\RR^d} |\mathcal{F}\Gamma_d(r)(\xi)|^2\mu(d\xi)dr&=& C\int_0^{\varepsilon} \int_{\RR^d}\frac{\text{sin}^2( r|\xi|)}{ |\xi|^2}
    ( |\xi|^2+1  ) ^{-\frac {\alpha}2} d\xi dr \\
&=& C\int_0^{\varepsilon} \int_{\RR^d}\frac{\text{sin}^2 |{\xi}|}{ |{\xi}|^2}\frac{r^{\alpha+2-d}}{(|{\xi}|^2+ r^2)^{ \frac {\alpha}2}}d{\xi} dr \\
&\geq& C\int_0^{\varepsilon} r^{\alpha+2-d} dr \int_{\RR^d}\frac{\text{sin}^2  |{\xi}| }{ |{\xi}|^2}\frac{1}{(|{\xi}|^2+1)^{\frac {\alpha} 2}}d{\xi} \\
&=& C\varepsilon^{\alpha+3-d}\,.
\end{eqnarray*}
Therefore, condition {(\bf H3)} is satisfied for $\eta=\alpha+3-d$.
To show (${\bf H 4}$)  as in the case of the heat equation we use the fact that  for any $x \in \RR^d$, $f(x)\leq C |x|^{-d+\alpha}$. Therefore, proceeding as in the case of the Riesz kernel with $\beta= d-\alpha$ we obtain that conditions \eref{hyp 4-1}, \eref{hyp 4-2} and  \eref{hyp 4-3} in (${\bf H 4}$)  hold, with $\eta_1 =  \alpha+3-d +\kappa_1$, $\eta_2=3 $ and  $\eta_3 = \alpha+3-d +\kappa_2$, respectively.

\medskip
\noindent
(C) \textit {Fractional  kernel}.
 Let $f(x)=\prod_{j=1}^d |x_j|^{2H_j-2}$, $ \frac{1}{2} < H_j < 1$ for $1\leq j \leq d$ such that $\sum_{j=1}^d H_j > d-1$.
Hypothesis ({\bf H1}) is verified by direct calculation.   By Section 3 in \cite{NuaE},  ({\bf H2}) holds when $d=1$   with $\kappa_1, \kappa_2 \in (0,H_1)$ and when $d=2$, it is satisfied for $\kappa_1, \kappa_2 \in (0, H_1+H_2-1)$.
  By Theorem 6.1 in \cite{HHN}, when $d=3$  ({\bf H2})  is satisfied with  $\kappa_1, \kappa_2 \in (0, \min(H_1+H_2+H_3-2, H_1-\frac{1}{2}, H_2-\frac{1}{2}, H_3-\frac{1}{2}))$. For ({\bf H3}), direct calculation and the change of variable $t\xi \rightarrow \xi$ yields
\begin{eqnarray*}
  \int_0^{\varepsilon} \int_{\RR^d} |\mathcal{F}\Gamma_d(t)(\xi)|^2\mu(d\xi)dt
&=&C \int_0^{\varepsilon} \int_{\RR^d}\frac{(\sin(t|\xi|))^2}{|\xi|^2}\prod_{j=1}^d |\xi_j|^{1-2H_j} d\xi dt\\
&=&C \int_0^{\varepsilon} t^{2\sum_{j=1}^d H_j-2d +2} dt \int_{\RR^d} \frac{(\sin (|\xi|))^2}{|\xi|^2}\prod_{j=1}^d |\xi_j|^{1-2H_j} d\xi\\
&=&C \varepsilon^{2\sum_{j=1}^d H_j-2d+3}\,,
\end{eqnarray*}
so ({\bf H3}) is satisfied with $\eta= 2\sum_{j=1}^d H_j-2d+3$.  For ({\bf H4}), we will check \eref{hyp 4-2} and \eref{hyp 4-3} first.  For \eref{hyp 4-2}, proceeding as before,
\begin{eqnarray*}
\int_0^{\varepsilon} r^{\kappa_1}\|\Gamma_d(r)\|^2_{\mathcal{H}} dr
&=& \int_0^{\varepsilon} \int_{\RR^d} r^{\kappa_1} \frac{(\sin(r |\xi|))^2}{|\xi|^2} \prod_{j=1}^d |\xi_j|^{1-2H_j}d\xi dr \\
&=&C \int_0^{\varepsilon} r^{\kappa_1+2\sum_{j=1}^d H_j -2d+2}dr = C \varepsilon^{\kappa_1+2\sum_{j=1}^d H_j -2d+3}\,,
\end{eqnarray*}
so \eref{hyp 4-2} is satisfied with $\eta_1=\kappa_1+2\sum_{j=1}^d H_j -2d+3$, which is strictly greater than $\eta$.  For \eref{hyp 4-3}, noting that the support of $\Gamma_d(r,*)$ is contained in the ball centered at the origin with radius $r$, we get
\begin{eqnarray*}
 \int_{0}^{\varepsilon} \langle |*|^{\kappa_2} \Gamma_d(r,*), \Gamma_d(r,\tilde{w}+*)\rangle_{\mathcal{H}}dr
&\leq&\int_{0}^{\varepsilon} r^{\kappa_2} \langle  \Gamma_d(r,*), \Gamma_d(r,\tilde{w}+*)\rangle_{\mathcal{H}}dr\\
&\leq& \int_{0}^{\varepsilon} r^{\kappa_2}  \int_{\RR^d} \frac{(\sin (r|\xi|))^2}{|\xi|^2} \prod _{j=1}^d |\xi_j|^{1-2H_j}d\xi dr\\
& =& C \varepsilon^{\kappa_2+2\sum_{j=1}^d H_j -2d+3}\,,
\end{eqnarray*}
so \eref{hyp 4-3} is satisfied with $\eta_3 = \kappa_2+2\sum_{j=1}^d H_j -2d+3$, which is strictly greater than $\eta$. For \eref{hyp 4-1}, we need to treat the cases $d=1,2,3$ separately. When $d=1$, fix $w \neq 0$. We have
\begin{eqnarray*}
  \int_0^{\varepsilon}\langle \Gamma_1(r,*), \Gamma_1(r,w+*)\rangle_{\mathcal{H}} dr
 = \frac{1}{4}\int_0^{\varepsilon} \int_{\mathbb{R}} \int_{\mathbb{R}} {\bf 1}_{\{ |x|< r\}} |x-y|^{2H_1-2} {\bf 1}_{\{ |y+w|< r\}} dy dx dr\,.
\end{eqnarray*}
When $\varepsilon $ is small enough, we  need to  have $|x-y| \geq C$ for some positive constant $C$  for the above  integrand to be  non zero. Hence,
\begin{eqnarray*}
 \int_0^{\varepsilon}\langle \Gamma_1(r,*), \Gamma_1(r,w+*)\rangle_{\mathcal{H}} dr
&\leq&C \int_0^{\varepsilon} \int_{\mathbb{R}} \int_{\mathbb{R}} {\bf 1}_{\{ |x|< r\}} {\bf 1}_{\{ |y+w|< r\}} dy dx dr\\
&=&C \int_0^{\varepsilon} r^2 dr = C \varepsilon^3\,,
\end{eqnarray*}
and when $d=1$ \eref{hyp 4-1} is satisfied with $\eta_2=3$, which is strictly greater than $\eta$.

When $d=2$, fix a nonzero point $w=(w_1,w_2) $. Without loss of generality, we may assume $w_1$ is not zero. We have
\begin{eqnarray*}
&&\int_0^{\varepsilon}\langle \Gamma_2(r,*), \Gamma_2(r,w+*)\rangle_{\mathcal{H}} dr\\
&=&\frac{1}{4\pi^2} \int_0^{\varepsilon} \int_{|x|< r} \int_{|y+w|< r} \frac{1}{\sqrt{r^2-|x|^2}} |x_1-y_1|^{2H_1-2} |x_2-y_2|^{2H_2-2} \frac{1}{\sqrt{r^2-|y+w|^2}} dx dy dr\,.
\end{eqnarray*}
Again, if $\varepsilon$ is small enough, we  must have $|x_1-y_1|> C$ for some positive constant $C$ for the above integral to be non zero.  Hence,   using the Fourier transform we obtain
\begin{eqnarray*}
&&\int_0^{\varepsilon}\langle \Gamma_2(r,*), \Gamma_2(r,w+*)\rangle_{\mathcal{H}} dr\\
&\leq& C \int_0^{\varepsilon} \int_{|x|< r} \int_{|y+w|< r} \frac{1}{\sqrt{r^2-|x|^2}}|x_2-y_2|^{2H_2-2} \frac{1}{\sqrt{r^2-|y+w|^2}} dx dy dr\\
&=&C \lim_{\delta \to 0} \int_0^{\varepsilon} \int_{|x|< r} \int_{|y+w|< r} \frac{1}{\sqrt{r^2-|x|^2}} e^{-\frac{\delta}{2}|x_1-y_1|^2}|x_2-y_2|^{2H_2-2} \frac{1}{\sqrt{r^2-|y+w|^2}} dx dy dr\\
&=&C \lim_{\delta \to 0} \int_0^{\varepsilon} \int_{\RR^2} \frac{(\sin(r |\xi|))^2}{|\xi|^2} p_{\delta}(\xi_1) |\xi_2|^{1-2H_2} e^{-i w\cdot \xi}d\xi dr\\
&\leq&C  \lim_{\delta \to 0} \int_0^{\varepsilon} \int_{\RR^2} \frac{(\sin(r |\xi|))^2}{|\xi|^2} p_{\delta}(\xi_1) |\xi_2|^{1-2H_2} d\xi dr \\
&=& C \int_0^{\varepsilon} \int_{\mathbb{R}} \frac{(\sin(r |\xi_2|))^2}{|\xi_2|^2} |\xi_2|^{1-2H_2} d\xi_2 dr =C \varepsilon^{2H_2+1}\,,
\end{eqnarray*}
so \eref{hyp 4-1} is satisfied with $\eta_2=\min(2H_1+1, 2H_2+1)$, which is strictly greater than $\eta$.

When $d=3$, fix a nonzero $w=(w_1,w_2,w_3) \in \RR^3$, without loss of generality, we may assume that $w_1 \neq 0$. We have
\begin{eqnarray*}
&&\int_0^{\varepsilon}\langle \Gamma_3(r,*), \Gamma_3(r,w+*)\rangle_{\mathcal{H}} dr\\
&=& \int_0^{\varepsilon} \int_{\RR^3} \int_{\RR^3} \Gamma_3(r,dx) \Gamma_3(r,w+dy) \prod_{j=1}^3|x_j-y_j|^{2H_j-2} dr\,.
\end{eqnarray*}
Again, when $\varepsilon$ is small enough, to make $x$ and $w+y$ in the support of the measure $\Gamma_3(r)$, we must have $|x_1-y_1|> C$ for some positive constant $C$. So
\begin{eqnarray*}
&&\int_0^{\varepsilon}\langle \Gamma_3(r,*), \Gamma_3(r,w+*)\rangle_{\mathcal{H}} dr\\
&\leq&C  \int_0^{\varepsilon} \int_{\RR^3} \int_{\RR^3} \Gamma_3(r,dx) \Gamma_3(r,w+dy) \prod_{j=2}^3|x_j-y_j|^{2H_j-2} dr\\
&=&C \lim_{\delta \to 0} \int_0^{\varepsilon} \int_{\RR^3} \int_{\RR^3} \Gamma_3(r,dx)\Gamma_3(r,w+dy)e^{-\frac{\delta}{2}|x_1-y_1|^2} \prod_{j=2}^3|x_j-y_j|^{2H_j-2} dr\\
&\leq&C \lim_{\delta \to 0}\int_0^{\varepsilon} \int_{\RR^3} \frac{(\sin r |\xi|)^2}{|\xi|^2}p_{\delta}(\xi_1)\prod_{j=2}^3 |\xi_j|^{1-2H_j} d\xi dr\\
&=&C \int_0^{\varepsilon} \int_{\RR^2} \frac{(\sin (r |(\xi_2,\xi_3)|))^2}{(\xi_2^2+\xi_3^2)}\prod_{j=2}^3|\xi_j|^{1-2H_j} d\xi_2d\xi_3 dr\\
&=&C \int_0^{\varepsilon} r^{2(H_2+H_3)-2} dr = C \varepsilon^{2(H_2+H_3)-1}\,,
\end{eqnarray*}
and \eref{hyp 4-1} is satisfied with $\eta_2=\min(2(H_2+H_3)-1, 2(H_1+H_3)-1,2(H_1+H_2)-1)$, which is strictly greater than $\eta$.


\begin{thebibliography}{99}

\bibitem{BP} Bally, V., Pardoux, E.: Malliavin calculus for white noise driven parabolic SPDEs. {\it Potential Anal.} 9 (1998), no. 1, 27-64.


\bibitem {CS} Chaleyat-Maurel, M., Sanz-Sol\'e, M.: Positivity of the density for the stochastic wave equation in two spatial dimensions. {\it ESAIM Probab. Stat.} 7 (2003), 89-114.

\bibitem{Da} Dalang,R.: Extending martingale measure stochastic integral with applications to spatially homogeneous s.p.d.e's. {\it Electron. J. Probab.} {\bf 4} (1999), 1-29.

\bibitem{DQ}Dalang, R., Quer-Sardanyons, L.: Stochastic integrals for spde's: a comparison. {\it Expo. Math. }29 (2011), no. 1, 6-109.

\bibitem{Gra} Grafakos, L.: Modern Fourier analysis.
Second edition. Graduate Texts in Mathematics, 250. {\it Springer, New York}, 2009. xvi+504 pp.

\bibitem{Fed} Federer, H.: Geometric measure theory.
Die Grundlehren der mathematischen Wissenschaften, Band 153 {\it Springer-Verlag New York Inc., New York} 1969 xiv+676 pp.

\bibitem{HHN} Hu, Y., Huang, J., Nualart, D.: On H\"older continuity of the solution of stochastic wave equations in dimension three. {\it Stoch. Partial Differ. Equ. Anal. Comput.} 2 (2014), no. 3, 353-407.


\bibitem{KX}Khoshnevisan, D., Xiao, Y., Harmonic analysis of additive Lvy processes. {\it Probab. Theory Related Fields} 145 (2009), no. 3-4, 459-515.


\bibitem{MNQ} Marinelli, C., Nualart, E., Quer-Sardanyons, L.: Existence and regularity of the density for solutions to semilinear dissipative parabolic SPDEs.
    {\it Potential Anal.} 39 (2013), no. 3, 287-311.



\bibitem{MS} Millet, A., Sanz-Sol\'e, M.: A stochastic wave equation in two space dimension: smoothness of the law. {\it Ann. Probab.} 27 (1999), no. 2,
    803-844.


\bibitem {NuaE} Nualart, E.: On the density of systems of non-linear spatially homogeneous SPDEs. {\it Stochastics} 85 (2013), no. 1, 48-70.

\bibitem{NuaD} Nualart, D.: {\it The Malliavin calculus and related topics.} Second edition. Probability and its Applications (New York). Springer-Verlag,
    Berlin, 2006. xiv+382 pp.


\bibitem{NQ} Nualart, D., Quer-Sardanyons, L.: Existence and smoothness of the density for spatially homogeneous SPDEs. {\it Potential Anal.} 27 (2007), no. 3,
    281-299.

\bibitem{QS1} Quer-Sardanyons, L., Sanz-Sol\'e, M.: A stochastic wave equation in dimension 3: smoothness of the law. {\it Bernoulli }10 (2004), no. 1, 165-186.


\bibitem{QS2} Quer-Sardanyons, L., Sanz-Sol\'e, M.: Absolute continuity of the law of the solution to the 3-dimensional stochastic wave equation. {\it J. Funct.
    Anal}. 206 (2004), no. 1, 1-32.

\bibitem{San} Sanz-Sol\'e, M.: Malliavin calculus.
With applications to stochastic partial differential equations. Fundamental Sciences. { \it EPFL Press, Lausanne; distributed by CRC Press, Boca Raton, FL}, 2005. viii+162 pp.

\bibitem{SS} Sanz-Sol\'e, M., Sarr\'a, M.: H\"older continuity for the stochastic heat equation with spatially correlated noise. {\it Seminar on Stochastic Analysis, Random Fields and Applications, III (Ascona, 1999)}, 259-268,
Progr. Probab., 52, {\it Birkhuser, Basel, 2002}.

\bibitem{Wal} Walsh, J.: An introduction to stochastic partial differential equations. {\it \'Ecole d'\'et\'e de probabilit\'es de Saint-Flour, XIVÑ1984}, 265-439,
Lecture Notes in Math., 1180, Springer, Berlin, 1986.

\end{thebibliography}
\end{document}